\documentclass[12pt,a4paper, reqno]{amsart}
\newcommand{\bt}{\begin{Theorem}}
\newcommand{\et}{\end{Theorem}}
\newcommand{\bi}{\begin{itemize}}
\newcommand{\ei}{\end{itemize}}
\newcommand{\bea}{\begin{eqnarray}}
\newcommand{\eea}{\end{eqnarray}}
\newtheorem{Definition}{Definition}[section]
\newtheorem{Theorem}[Definition]{Theorem}
\newtheorem{Lemma}[Definition]{Lemma}

\newtheorem{Proposition}[Definition]{Proposition}
\newtheorem{Corollary}[Definition]{Corollary}

\newtheorem{Remark}[Definition]{Remark}
\newcommand{\be}{\begin{equation}}
\newcommand{\ee}{\end{equation}}
\newcommand{\cD}{\mathcal{D}}

\newtheorem{exam}{Example}

\newtheorem{theorem}[exam]{Theorem}

\newcommand{\Dl}{\Delta}

\begin{document}

\title [$CR$-foliations]{Propagation of boundary 
$CR$ foliations and Morera type theorems for manifolds
with attached analytic discs}
\author { Mark L. ~Agranovsky}
\thanks { This work was partially
supported by Israel Scientific Foundation,
grant No. 279/02-01.}

\address{ M.L. Agranovsky, Department of
Mathematics, Bar-Ilan University, 52900, Ramat Gan,
Israel}
\email{agranovs@macs.biu.ac.il}



\vskip.25in

\setcounter{equation}{0}

\maketitle
\begin{abstract}

We prove that homologically nontrivial generic smooth $(2n-1)$-parameter
families of analytic discs in
$\mathbb C^n, \ n \le 2,$ attached by their boundaries
to a $CR-$ manifold $\Omega,$ test $CR$-functions
in the following sense: if a smooth function on $\Omega$ analytically extends
into any analytic discs from the family,
then the function  satisfies tangential $CR-$ equations on $\Omega.$

In particular, we give an answer (Theorem 1) to the
following long standing open question, so called strip-problem, 
earlier solved only for special families (mainly for circles):
given a smooth one-parameter family of Jordan curves in the plane 
and a function $f$ admitting holomorphic extension inside each curve, 
must $f$ be holomorphic on the union of the curves? 
We prove, for real-analytic functions and arbitrary generic real-analytic 
families of curves, that the answer is ``yes'', if no point is surrounded 
by all curves from the family. The latter condition is essential.
We generalize this result to characterization of complex curves
in $\mathbb C^2$ as real 2-manifolds  admitting nontrivial
families of attached analytic discs (Theorem 4.)

The main result implies fairly general Morera type 
characterization of $CR$-functions on hypersurfaces in $\mathbb C^2$
in terms of holomorphic extensions into three-parameter families of
attached analytic discs (Theorem 2).
One of the applications is confirming, in real-analytic category, 
the Globevnik-Stout 
conjecture (Theorem 3) on boundary values of holomorphic functions. 
It is proved that a smooth function on the boundary of a smooth 
strictly convex domain in $\mathbb C^n$ extends holomorphically 
inside the domain if it extends holomorphically into     
complex lines tangent to a given strictly convex subdomain.

The proofs are based on a universal approach, namely, on the
reduction to a problem of propagation, from the boundary to the interior,
of degeneracy of $CR$-foliations 
of solid torus type manifolds (Theorem 2.2).

\end{abstract}

\vskip 0.5cm

\section{Formulation of the problem, the main results and comments.}\label{1}

\noindent
{\bf 1.1. Formulations of the problems and definitions.} 

\medskip
The results of this article are related to the
following general problem:

\bigskip
\noindent
{\bf Morera problem  for $CR$-functions:} 
{\it Let $\Omega \subset \mathbb C^n$ be a $CR-$ manifold of real
dimension $k,$
and let $\cD$ be a family of analytic discs in $\mathbb C^n$,
such that the boundaries
$\partial D, \ D \in \cD,$ cover $\Omega.$ Let $f$ be a continuous or smooth
function on $\Omega$ such that for each $D \in \cD$
the restriction $f\vert_{\partial D}$ admits analytic extension in $D.$
When does this imply that $f$ is a $CR-$ function in $\Omega?$}

Recall that the differentiable
manifold $\Omega$ is called $CR-$ manifold if the
dimension of maximal complex subspaces $T^{\mathbb C}_p(\Omega)$ of the
tangent spaces are the same  for all points $p \in \Omega$. A smooth function
$f$ on $\Omega$ is called $CR-$ function if $\overline {\partial}_b f=0$
for any tangential Cauchy-Riemann operator $\overline {\partial}_b$ in the
complex tangent space.

The main result of this article is Theorem 2.2 which implies an
answer the above question
in the cases $n=1, \ k=2$ and $n=2, \ k=3.$ 
We consider smooth regular $(2n-1)-$ parameter
families of discs and obtain conditions for the families to detect $CR-$
functions.

The condition for the family $\cD$ involves the homology  or
relative homology groups, depending on whether the parameterizing family is
closed or not.
For the planar case, $n=1, \ k=2,$
the condition for the family is most simple:
the closures of the discs must have no common point.

In the case $n=1, \ k=2$ our main result implies 
a general solution, in real-analytic category, 
of the following 
problem which was known for a long time and was
solved till now only for special families of curves:
\footnote{ We require real-analyticity to guarantee 
nice structure of zero sets
of certain functions exploited in our constructions
(see Section 3.7). This assumption seems redundant and likely might be
reduced to conditions of differentiablity by more refined analytic arguments,
in particular, by using the technique of currents (see the concluding remarks
at the end of the article).}

\medskip
\noindent
{\bf Strip-problem.} {\it Let $\gamma_t$ be  a continuous (smooth)
one-parameter family of Jordan curves
in the complex plane. Let $f$ be a continuous (smooth) function
on the union $\Omega= \cup \gamma_t$. Suppose that for each $t$
the restriction
$f$ admits analytic extension inside the curve $\gamma_t.$
When does this imply that $f$ is analytic in $\Omega?$}

This question naturally arose in relation with the works \cite{AV},\cite{Z1}, \cite {G1}
and, to the knowlegde of the author, first was explicitly formulated by
Josip Globevnik in his talk \cite{JG}.

Note that the strip-problem is, in a sense, a degenerate case
when the attached analytic discs lie in the same complex plane to which
the manifold $\Omega$ belongs. 

\medskip
For the dimensions $n=2, \ k=3,$ i.e. for real hypersurfaces $\Omega$
in $\mathbb C^2$, 
we prove quite general Morera type theorem for $CR-$ functions on $\Omega.$
Using this theorem we
prove, in real-analytic category, the following
conjecture formulated by J.Globevnik and
E. L. Stout in the article \cite{GS}  (independently this question was posed by the author
in the talk \cite{MLA}):

\medskip
\noindent
{\bf  Globevnik-Stout conjecture  (one-dimensional extension property).}
{\it Let $D$ be a strictly convex
bounded domain in $\mathbb C^n$ and $S \subset D$
be a (smooth convex) closed hypersurface, compactly belonging to
$D.$
Suppose that a continuous (smooth) function $f$ on $\partial D$ possesses
the property: for any complex line $L$ , tangent to $S,$ the restriction
$f\vert L \cap \partial D$ analytically extends to $L \cap D$. Then
$f$ is the boundary value  of a holomorphic function in $D.$}

\bigskip
\noindent
{\bf 1.2. Formulations of the main results.} 

\medskip
First of all, let us describe in a more precise way 
the families of curves and analytic
discs we are going to deal with.

Let $M$ be a compact connected smooth oriented $(2n-1)-$manifold with
boundary $\partial M$
which, in particular, may be empty. For instance, for $n=1$ the manifold
$M$ is a smooth curve, topologically equivalent either to the circle
$S^1$ or to the closed segment $[0,1].$
In this case,we will think of $M$ as of the unit circle $M=S^1$ in
the complex plane, or the segment $M=[0,1].$

By {\it analytic disc} in the complex space $\mathbb C^n$ we understand
a holomorphic diffeomorphic embedding
$$g:\Delta \mapsto \mathbb C^n,$$
of the unit disc $\Delta$  in the complex plane. The mapping
$g$ is assumed smooth up to $\partial \Delta$ and
$$g^{\prime} (\zeta) \neq 0, \forall \zeta \in \overline \Delta.$$
Sometimes we will use the term ``analytic disc'' for the image $D=g(\Delta).$

Given a real manifold $\Omega \subset \mathbb C^n$, the analytic disc $D$
is called {\it attached} to $\Omega$ if its boundary
$\partial D=g(\partial \Delta) \subset \Omega.$ It may happen that the entire
analytic disc $D$ or its portion belongs to $\Omega.$

By the {\it $C^r$-family of analytic discs} parameterized by the manifold
$M, dim M=k <2n,$ we understand the family  $D_t=g_t(\Delta),$
\ where $g_t, \ t \in M,$ is a family of holomorphic embeddings 
(or even immersions)
of the disc $\Delta,$ smoothly parametrized by the points in $M.$
Here we assume that
$$G(\zeta,t)=g_t(\zeta), \  (\zeta,t) \in \Delta \times M,$$
belongs to the class $C^r(\overline \Delta \times M) \ \ r \ge 2.$ 
If the manifold $M$ is real-analytic and the mapping
$G(\zeta,t)$ is real analytic in the closed domain
$\overline \Delta \times M$ then we say that the mapping $G$
parameterizes a {\it real-analytic family} of analytic discs.
In most considerations we will need only differentiability of the manifold
$M$ and of the parameterization mapping  $G.$


Denote $\Omega$ the set covered by the boundaries of the analytic discs:
$$\Omega=G(\partial \Delta \times M)=
\cup_{t \in M} \gamma_t.$$ Regularity assumptions for $G$,
which we will discuss below,
provide that $\Omega$ is a smooth manifold.
We also denote
$$\hat \Omega=G(\overline \Delta \times M)= \cup_{t \in M} \overline{D_t}.$$

Throughout the article we will use the notations:
$$\Sigma=\Delta \times M, \ \ b\Sigma=\partial \Delta \times M, \ \ 
\Sigma_0=\Delta \times \partial M, \ \ 
\partial \Sigma= b\Sigma \cup \Sigma_0.$$

In this article we consider the  following two cases:

\medskip
\noindent
{\bf The case A, n=1.}

In this case the analytic discs $D_t$ are
Jordan domains, bounded by Jordan curves $\gamma_t=\partial D_t,$
in the complex plane and $\Omega \subset \mathbb C.$ The manifold $M$ has
the dimension $2n-1=1$ and therefore $M$ is topologically
either a circle $S^1$ or a closed segment $[0,1] \subset \mathbb R.$

We will assume the regularity condition. The family $\{D_t\}_{t \in M}$
of Jordan domains will be called
{\it regular} if the parameterizing mapping
$G$ has the minimal degeneracy:

\noindent
1. $rank \ dG (p) = 2, \ \forall p \in \Sigma.$

\noindent
2. $rank \ dG\vert_{b\Sigma}(p)=2,
\ \forall p \in b \Sigma \setminus  Crit(G),$
where $Crit(G) \subset b \Sigma$
is the one-dimensional critical smooth manifold of $G.$

\noindent
3. $rank \ dG\vert_{b \Sigma}(p)=1, \ \forall p
\in Crit(G),$ and $G(Crit(G)) \subset \partial \Omega.$

Under these conditions,
the set $\Omega$ is a closed domain,
its interior points are regular values for the restriction of
$G$ to the boundary manifold $b \Sigma$ and
the critical values are located on the boundary.
The preimage $G^{-1}(\partial \Omega)\cap b\Sigma$ 
contains
the critical curve $Crit(G)$ and the mapping $G$ is regular  out
of the critical curve.

\medskip
\noindent
{\bf The case B, n=2.}

In this case $\Omega \subset \mathbb C^2.$
We are interested in the case when
$\Omega$ is a smooth  real hypersurface in $\mathbb C^2 \cong
\mathbb R^4.$
Then by {\it regularity} of the family $D_t$ is understood the
following conditions:

\noindent
1. $rank \ dG (p) =4, \ \forall p \in \Sigma= \Delta \times M.$

\noindent
2. $rank \ dG\vert_{b \Sigma}(p)=3,
\ \forall p \in b \Sigma \setminus  Crit(G),$
where $Crit(G) \subset b \Sigma$ is the two-dimensional smooth
critical manifold.

\noindent
3. $rank \ dG\vert_{b \Sigma}(p)=2, \ \forall p
\in Crit (G), $ and $G(Crit(G)) \subset \partial \Omega.$

Thus, the non-boundary points of the manifold $\Omega$
are regular values for the mapping $G : b \Sigma \mapsto
\Omega,$ while the critical values lie in $\partial \Omega$.

\bigskip
Now we can formulate the main results.

\noindent
{\bf The results for the case  A ($\Omega$ is a domain in $\mathbb C$).}

In the case $n=1$ the parametrizing manifold $M$ is diffeomorphic to
either unit circle $S^1$ or the segment $[0,1].$ 
We assume that $M$ is just one of these curves.

\begin{theorem} Let $\{D_t\}_ {t \in M}, \ \ M=S^1$ or $M=[0,1],$
be a real-analytic regular closed family of Jordan domains in
$\mathbb C.$

Assume that the closures $\overline D_t$ have no common point:

\noindent
(a)  $\cap_{t \in M} \overline {D_t} = \emptyset .$

Let $f$ be a real-analytic function in $\overline \Omega$ and assume that
$f$ satisfies the property:

\noindent
(*)   for each $t \in M$ the restriction $f\vert_{\partial D_t}$
admits holomorphic extension in  $D_t.$

Then $f$ is holomorphic in the interior of
$\Omega=\cup_{t \in M} \partial D_t$ and extends holomorphically
to $\hat \Omega=\cup_{t \in M} D_t.$

The condition
(a) cannot be omitted.
\end{theorem}

\medskip
\noindent
{\bf The results for the case B ($\Omega$ is a hypersurface in $\mathbb C^2$).}

The condition (a) in general form is the following. 

\begin{Definition} Let $\partial M=\emptyset.$
We say that the family of the analytic discs
$\cD= \{D_t\}_{t \in M},$ parametrized by the mapping
$G(\zeta,t), \ (\zeta,t) \in \Sigma=\Delta \times M,$
is homologically nontrival if the induced mapping of the relative
homology groups
$$G: H_k (M,\partial M;\mathbb R) \cong 
H_k(\overline \Sigma, \overline \Sigma_0; \mathbb R) \mapsto 
H_k(G(\overline \Sigma), G(\overline \Sigma_0); \mathbb R),$$
$k =\dim \ M$, is not trivial, $G_* \neq 0.$
\end{Definition}

Let us give a geometric interpretation of the condition of homological 
nontriviality.
For the case $\partial M =\emptyset$ it reads as follows:
no $d$-cycle $c \subset \hat \Omega=\cup_{t \in M} \overline D_t, \ 
d =\dim M,$ meeting any closed analytic dics $\overline D_t, t \in M,$
is the boundary
$c = \partial c^{\prime}$ of a $(d+1)$- cycle
$c^{\prime} \subset \hat \Omega.$

In the case $\partial M \neq \emptyset$ (``non closed family"') 
the condition means the following: if a
$d$-cycle $c \subset \hat \Omega$ intersects
each closed disc $\overline D_t, \ t \in M,$ then $c$ is relatively
homologically nontrivial, i.e. for no $d$-cycle $c_1 \subset
\cup_{t \in \partial M} \overline D_t$ the union $c \cup c_1$
is a boundary of a $(d+1)$- cycle $c^{\prime} \subset \hat \Omega.$

\begin{theorem} Let $\cD= \{D_t\}_{ t \in M}$ be a
real-analytic regular family
of analytic
discs in $\mathbb C^2,$ parametrized by a real connected
3-dimensional real analytic compact manifold
$M$ with boundary (possibly empty).  
Denote $\Omega=\cup_{t \in M} \partial D_t-$
a real-analytic  real 3-dimensional closed manifold in $\mathbb C^2.$
Assume that

\noindent
(a) the family $\cD$ is  homologically nontrivial.

Let $f$ be a real-analytic function on $\Omega$  such that

\noindent
(*) for each $t \in M$ the restriction $f\vert_ {\partial D_t}$
admits analytic extension in the analytic disc $D_t.$

Then $f$ is a $CR-$ function on $\Omega.$
\end{theorem}
%
%
%

Recall that $\hat \Omega$ is defined as 
the union of all closed analytic discs $\overline{ D_t}.$ 

\noindent
{\bf Remark.}
The condition (*) is equivalent to the following:
$$\int_{\partial D_t} f \omega=0,$$
for any holomorphic differential 1-form $\omega$ in $\mathbb C^2.$ 
\noindent
2. In Section 2.2 we will show that in the planar case 
the homological nontriviality is equivalent to the condition (a) in 
Theorem 1 that the closed domains $\overline D_t$ have no common point.

\medskip
Since for $n>1$ the condition of holomorphic extendibility is local
(differential) as opposite to the case $n=1$ when it is a global
(integral) condition,
the boundary values of holomorphic functions can be checked by complex
two-dimensional sections and our result for  $\mathbb C^2$ imply
corresponding characterization of boundary values of holomorphic functions of $n$
variables for $n \ge 2.$
So, Theorem 3 implies answer to the
question of of Globevnik and Stout \cite{GS} in real-analytic category:

\begin{theorem}
Let $D \subset \mathbb C^n$
be a bounded strictly convex domain with real-analytic boundary
and $S \subset D$ be a real-analytic strictly convex 
hypersurface. 
Let $f$ be a real-analytic function on $\partial D$ such that
for any complex line $L$ tangent to $S$  the restriction
$f\vert L \cap \partial D$ extends analytically in $L \cap D.$
Then $f$ is the boundary value of a function holomorphic in $D$ and
continuous
in $\overline D.$
\end{theorem}

\begin{proof} If $n=2,$ then Theorem 3 is  a particular case
of Theorem 2, for a special family of
sections by complex lines. 
In this case the parametrizing manifold
is $M=S.$ The analytic disc $D_t, t \in S,$ is defined as the intersection 
$$D_t=L_t \cap D$$
of the domain $D$ with the complex line $L_t$ tangent ot $S$ at the point $t.$

Then $$\hat \Omega= \cup_{ t \in S} \overline {D_t}= \overline D \setminus D^{\prime},$$
where $D^{\prime}$ is the domain bounded by $S.$ 
If a 3-cycle $c$ intersects each
analytic disc $D_t$ then any 4-cycle 
bounded by $c$ must contain $D^{\prime}$ and hence is not contained in $\hat \Omega$. Therefore
$c$ is not homologous to zero in $\hat \Omega$ and the condition 
(a) of Theorem 2 holds. Therefore $f$ is $CR$ function on $\partial D$ and hence is the boundary value of
a holomorphic function in $D.$   

If $n>2,$ then consider complex 2-planes $\Pi$ intersecting
$S.$ For almost all $\Pi$ the intersection $S \cap \Pi$ is
a real-analytic hypersurface, contained in the intersection
 $$D_{\Pi}= D \cap \Pi$$ which
can be regarded as a strictly convex domain in $\mathbb C^2.$
The surface $S \cap \Pi$ is strictly convex.

Thus, we are in position of Theorem 2 applied to the domain
$D_{\Pi} \subset \mathbb C^2$
and to the family of complex lines $L \subset \Pi,$ tangent to the 3-surface
$S_{\Pi}.$
By Theorem 3 for $n=2$ the function $f$ is annihilated by
any tangential differential operator $\overline X, \ X \in T_p^{\mathbb C}(\partial D \cap \Pi), \  p \in \partial D \cap \Pi.$
Due to the large supply of the complex
2-sections $\Pi,$ we obtain $\overline X f =0$ for any point $p \in \partial D$ and any vector $X \in T^{\mathbb C}(\partial D)$. 
Therefore $f$ is $CR$ function on $\partial D$ and hence  extends holomorphically in the domain $D.$
\end{proof}

\bigskip
\noindent
{\bf 1.3. Generalization of Theorem 1: characterization of 
complex curves in $\mathbb C^2.$}

\medskip
Theorem 1 can rephrased in terms of the graph $\Gamma_f$ of the function $f.$
Indeed, if $F_t$ is the analytic extension of $f$
from the boundary $\partial D_t$ into domain $D_t$ then the graph $\Gamma_{F_t}$ of $F_t$ over $D_t$ is an analytic disc in 
$\mathbb C^2$ attached to $\Gamma_f.$ Therefore Theorem 1 states that if
the real 2-manifold $\Gamma_f \subset \mathbb C^2$  admits
a nontrivial one-parameter family of attached analytic discs then $\Gamma_f$ is a complex manifold. 

This version of Theorem 1 can be generalized to compact manifolds which 
are not necessarily graphs.

For simplicity we formulate the result for periodic families of 
attached analytic discs:

\begin{theorem} Let $\Lambda \subset \mathbb C^2$ be a 
real-analytic 2-dimensional oriented compact manifold
with nonempty boundary, $\partial \Lambda  \neq \emptyset.$ Suppose that
$\Lambda$ admits one-parameter regular homologically nontrivial real-analytic 
family $D_t, \ t \in S^1,$ of attached analytic discs such that 
$\Lambda=\cup_{t \in S^1} \partial D_t.$
Then $ \Lambda $ 
is a 1-dimensional complex manifold in $\mathbb C^2$ 
and $D_t \subset \Lambda$ for all $t.$ 

\end{theorem} 

By regularity we understand here that $\Lambda \setminus
\partial \Lambda$ consists of  regular values of the parametrizing mapping
on $b \Sigma.$


\bigskip

\noindent
{\bf 1.4. History of the problems.} 

\medskip
\noindent
{\bf Strip-problem.} The name ``strip-problem'' is due to the 
typical shape of domains
swept out by one-parameter families of curves in the plane
(see,e.g.
Ehrenpreis' book \cite{E1}, p.575; for the strip-problem for more general PDE see  \cite{E1}, Ch.9.5,
and \cite{AN}). The analytic extendibility inside
a planar Jordan curve can be formulated it terms of a complex moments
condition, thus the question can be regarded
as a version of Morera theorem, in which the lowering by one  of the
number of parameters for the family  of the testing contours
is compensated by
the stronger condition of vanishing of {\it all} complex moments.
We also refer the reader for this and related problems to \cite{Z1},
\cite{Z2}.

In the paper \cite{AV}, by Val'sky and the author,
on Moebius-invariant function
algebras in the unit disc, a lemma was proved about testing of analyticity
by analytic extendibility into families of Jordan curves in the disc.
The families were assumed
invariant with respect to conformal automorphisms of the unit disc.
The method used was averaging of a function with respect to rotations
and applying the argument principle to the averaged function.

Globevnik \cite{G1}
observed that replacing the
averaging by computing the Fourier
coefficients in the polar coordinates leads to an analogous test of
analyticity for rotation-invariant families of curves.
In the articles \cite{G3},\cite{G4},\cite{G6} he made several
interesting observations on the phenomena.
In \cite{AS} the result of Globevnik \cite{G1} was generalized for
$U(n)-$ invariant families of boundaries of analytic discs in $\mathbb C^n,$
with using decompositions into spherical harmonics in $\mathbb C^n$.

The above results used tools of harmonic analysis and
therefore required  group invariance of the testing families.
However, even for (noncompact) group-invariant families of curves,
when no information about the growth of functions in the question is known and
Fourier analysis becomes inapplicable,
simple natural questions remained unanswered.

For instance, the following question
became a challenge: given a continuous or smooth function $f$ in the strip
$|Im z| < 1$, does the  analytic extendibility inside any
inscribed circle imply that
$f$ is holomorphic in the strip?

The first result beyond
harmonic analysis was obtained by Globevnik and the
author \cite{AG}. The problem was completely solved for arbitrary
one-parameter families of circles in the plane, though for functions
$f(x,y)$ which are rational (quotient of two polynomials) in $x,y.$
In spite of yet geometric restrictions for the curves (circles),
the approach in \cite{AG}  led to a new insight. The key point
was reformulating the problem, originally one-dimensional,
in $\mathbb C^2$ using the
embedding  $z \mapsto (z,\overline z)$ of the real 2-plane into
$\mathbb C^2$,  
The functions, along
with their analytic extensions, lift to the quadrics
$(z_1-a)(z_2-\overline a)=r^2$
in $\mathbb C^2$,
which are the complexifications of the circles $|z-a|=r.$
Then the proof in \cite{AG} is based on analysis of the
dynamics, in the parameter $t$, of the quadrics,
parameterized by $a=a(t)$ and $r=r(t)$,
with respect to the zero varieties of the
polynomials generating the rational function $f.$

Also, in \cite{AG} the case of real-analytic functions and arbitrary smooth
families of circles
was solved (independently , but in a special case, the same result
was obtained by Ehrenpreis \cite{E2}, with the help of Fourier analysis.)

The next significant progress was due to Tumanov \cite{T1}.
As in \cite{AG}, he also
started with the lifting the
problem into $\mathbb C^2$, but applied powerful tools of $CR-$ theory,
in particular, the edge of the wedge Ayrapetyan-Henkin's theorem, to prove
forced analytic extendibility of the lifted function $f$ to a larger domain.
Note, that in \cite{AG} such extension was provided automatically,
as rational functions in $z,\overline z$ always possess
meromorphic extensions inside any circle.

As the result, in \cite{T1} the
strip-problem was solved for
continuous functions, albeit for more restricted in the sense of \cite{AG}, 
families
of circles, namely those with constant radius and centers on an interval.
However, soon afterward, Tumanov \cite{T2}
got rid of the above restrictions and came up with a proof for the case
of continuous functions $f$ and arbitrary smooth families of circles.
Moreover, this proof, motivated by an argument of Hans Lewy,
was much simpler than that in \cite{T1}. Recently Globevnik \cite{G7}
generalized
the geometric construction from \cite{T1} and used the reasoning from
\cite{T2} to solve the strip-problem for special
families of non-circular Jordan curves which are
translates of a fixed axially symmetric Jordan curve,
along the line orthogonal to the symmetry axis.

In Theorem 1  of this article we give the solution for generic
families of general Jordan curves
with no restriction of geometric type.
Our approach rests on  a
reformulating of the original problem to the topological language,
namely, as a question about
$CR-$ extensions of coverings of the 2-dimensional
torus or cylinder, inside the solid torus or solid cylinder.

This reduction reveals
topological or, better to say, topology-analytical nature
of the problem, as well as the adequate tools for the solutions.
As result, it allows to get rid of geometric restrictions for the Jordan
domains in the question and solve the problem for general families.
Some ingredients of the analytic parts of the proofs are close to those
in the  article \cite{T2} by Tumanov.

\medskip
\noindent
{\bf One-dimensional extension property.}
We refer the reader for boundary Morera theorems to
the recent survey by Kytmanov and Myslivets \cite{KM},
and an extended  bibliography there. Here we will outline only
some results which are mostly related to our paper.

It was observed in \cite{AV} that boundary values of holomorphic functions
in the unit ball in $\mathbb C^n$ can be characterized by analytic 
extendibility
into sections of the  ball by complex lines. Stout \cite{St2} generalized this
result to arbitrary smooth domains, using complex Radon transform.
In \cite{AS1} the family of lines was
reduced to the set of complex lines passing through a fixed open domain.

Nagel and Rudin \cite{NR} proved one-dimensional property for the ball
in $\mathbb C^n$
and the family of complex lines tangent to a smaller concentric ball.
Globevnik \cite{G2,G5} reduced the families of lines in the question.
In \cite{ABC} an analogous result was obtained for lines tangent to a fixed 
orbit of the Heisenberg group acting on the complex ball.

The fundamental work by Globevnik and Stout \cite{GS}
contains many deep results on the subject. There the approach is mainly
based on the complex Radon transform in its various versions, in particular,
approximation by the complex plane waves. Tumanov \cite{T3} obtained
similar characterizations of $CR-$ functions on $CR-$ manifolds of higher
codimension.

Note, that of most interest are the families not
containing small analytic discs. By small discs we understand the 
families of discs that can be shrinked to boundary points.
In this case, the problem of testing $CR-$ functions is much easier, 
at least for
smooth functions, as the tangential
$CR-$ equation follow from  Stokes formula applied to the shrinking discs.

An example of families without small discs is the family of intersections
of a fixed domain with complex curves
(for instance, complex lines, as in \cite{NR}), tangent
to a fixed surface.
In \cite{GS} Globevnik and Stout
conjectured that the Nagel-Rudin theorem should be true
for two arbitrary  enclosed convex domains $D^{\prime} \subset D$,
when the family
of analytic discs, testing boundary values of holomorphic functions
on $\partial D,$ consists of sections of the domain $D$ by the complex lines tangent
to $\partial D^{\prime}.$

Dinh \cite{D} confirmed the conjecture for smooth functions and
under certain condition of ``strongly non real-analyticity'' of the above sections by the complex lines.

Recently, the conjecture of Globevnik and Stout was confirmed in
affirmative by Baracco, Tumanov and Zampieri \cite{BTZ}
for the family of extremal discs.
The extremal discs are geodesics
in Kobayashi metrics in the larger domain. The proof goes back to
the idea of the proof in \cite{T1} and hence the extremal discs are needed
for meromorphic lifts to tangent spaces, similarly to lifting
circles to complex quadrics in $C^2$ by means of $\overline z.$

In this article we prove (Theorem 2)
that in dimension $n=2$ and under assumptions of sufficient
smoothness (real-analyticity),
no geometric restrictions for analytic discs and for
character of the family are required, and
the one-dimensional extension property, at least for smooth functions,
is true for arbitrary generic family of attached analytic discs.
The essential condition is rather topological and requires
homological nontriviality of the family.

Since for $n>1$ characterization of (smooth)
boundary values of holomorphic functions are differential (local), as opposed
to $n=1$,
Theorem 2, albeit is stated for $n=2,$
lead to various boundary Morera theorems
in arbitrary dimensions. To apply Theorem 2  
the families of attached analytic discs must have  rich enough supply of
subfamilies filling up 2-dimensional complex submanifolds. Theorem 3 
is just an example of such kind of
application.

Theorems 1 and 2 are proved by an universal approach which we explain in the
following section. The remarkable 
papers of Alexander and Wermer \cite {AW} and
of Stout \cite {St1} about linking numbers and analytic functions
were a major influence in my discovering this approach.

\bigskip
\section{Reduction of the problems to $CR-$extensions
of coverings and foliations }\label{2}

\noindent
{\bf 2.1. Formulation of the equivalent result on extensions of foliations.}

\medskip
We will formulate a single theorem which includes, in equivalent form,
Theorems 1 and 2.

Let $n$ be an integer and $M$ be a compact connected  $C^r-$ smooth
oriented real $(2n-1)-$
manifold with the boundary $\partial M$, possibly empty.

As above, we denote
$$\Delta=\{z \in \mathbb C : |z|<1 \}.$$ For the unit circle in the plane
 we will use
both notations, $\partial \Delta$ and $S^1,$ depending on the context,
analytical or topological.
We also use all the notations, $\Sigma, \ b \Sigma, \Sigma_0, \ \partial \Sigma, \ \Omega, \ \hat \Omega,$ 
 from Section 1.
The real dimensions of the manifolds are
$$\dim_{\ \mathbb R} \ \Sigma=2n+1, \ \ \dim_{\ \mathbb R}  b\Sigma=2n.$$

The dimensions of the manifolds $\Omega=G(b \Sigma)$ and $\hat \Omega=
G(\overline \Sigma)$ 
depend on the rank of the mapping $G,$ which we are now going to specify. 

In order to combine both cases A and B, discussed in the previous section,
let us introduce the definition.

\begin{Definition} Let $k \le 2n$ be an integer. We say that a
smooth mapping
$$ G:\overline \Sigma \mapsto \mathbb C^n$$
is {\bf $(2n,k)-$  regular} if

\noindent
1. $\hat \Omega=G(\overline \Sigma)$ is a closed domain in $\mathbb C^n,$ \ \
$rank_{\ \mathbb R} \ dG(p)=2n$ for all $p \in \Sigma,$

\noindent
2. $\Omega=G(b \Sigma)$ is a smooth real $k$- dimensional
manifold in $\mathbb C^{n} \cong \mathbb R^{2n},$ with boundary \ \ and
$rank _{\ \mathbb R} \ dG\vert_{b\Sigma} (p)=k$ for all
$p \in b\Sigma \setminus Crit(G),$
where $Crit(G) \subset b \Sigma$ is $(k-1)$-dimensional
critical manifold, and $Crit(G) \subset G^{-1}(\partial \Omega),$

\noindent
3. $rank_{\ \mathbb R} \ dG\vert_{b \Sigma}(p)=k-1,$
for $p \in Crit(G).$

\end{Definition}

In this article we are concerned with the two cases

\bigskip
\noindent
{\bf (A)}: $n=1, \ k=2.$

In this case the image
$\hat \Omega=G(\overline \Sigma)$ has real dimension 2 and
is a closed domain in $\mathbb C$,
$G(b \Sigma)=\Omega$  is also 
a domain in $\mathbb C,$ with the smooth boundary
$\partial \Omega$ containing
the critical values of $G$. The mapping
$\zeta \mapsto G(\zeta,t)$ maps conformally the unit disc $\Delta$
onto the analytic disc $D_t \subset \mathbb C.$
These analytic discs are attached to $\Omega$ which means
that $\partial D_t \subset \Omega.$

The restriction $G\vert_{b \Sigma}$ of $G$
to the 2-dimensional boundary manifold
$b \Sigma $ is a finitely sheeted covering
over $\Omega \setminus \partial \Omega.$ The fibers correspond
to the curves 
$$\gamma_t=\partial D_t=G(\partial \Delta \times \{t\}),$$
passing through a fixed point in $\Omega.$

The set
$Crit(G) \subset G^{-1}(\partial \Omega)$ of critical points is a curve in
$b \Sigma.$
The set of critical values 
$G(Crit(G)) \subset \partial \Omega$ is the envelope of the family
of the curves $\gamma_t.$
The points from the envelope are the {\it sliding}
points. At these points the velocity vector (the motion direction)
of the family $\gamma_t$ is proportional to the tangent vector to $\gamma_t.$
Except the sliding points the boundary $\partial \Omega$ contains
subarcs of the curves $\gamma_t$ corresponding to
``boundary'' parameters $t \in \partial M.$  

The condition that the interior of $\Omega$ consists of regular
values of the mapping $G\vert_{b\Sigma}$ 
(does not contain sliding points) guarantees 
that this mapping is a foliation over the interior of $\Omega$. This condition is important in our considerations.

\bigskip
\noindent
{\bf (B)}: $n=2, \ k=3.$

In this case  $\hat\Omega$ has real
dimension $2n=4$, while $\Omega$ is a real
3-dimensional  submanifold of $\mathbb C^2,$ contained in
$\partial \hat\Omega.$ The analytic discs $D_t$ are attached to $\Omega.$
The union of their closures constitutes $\hat \Omega.$

The restriction of the parameterization mapping $G$
to the 4-dimensional manifold $b \Sigma$ is a foliation with the
3-dimensional base $\Omega \setminus \partial \Omega$ and the 1-dimensional fibers 
$G^{-1}(b), b \in \Omega \setminus \partial \Omega.$
The fibers correspond to the one-parameter family of the curves passing through
one point in $\Omega.$

\bigskip
Now we are ready to formulate the main
result, which is just a general form of Theorems 1 and 2.

\begin{Theorem}\label{st:thm}

Let $Q=(F,G) : \overline \Sigma \mapsto \mathbb C \times \mathbb C^n=
\mathbb C^{n+1}$
be a real-analytic $CR-$ mapping, that is $Q$ is holomorphic on each complex
fiber $\Delta \times \{t\}.$

Suppose that

\noindent
(i) The mapping $G$ is $(2n,k)-$ regular,

\noindent
(ii) The mapping $Q$ degenerates on the boundary
$b \Sigma$, meaning that
$$ F= f \circ G$$
holds on $b \Sigma$ for some smooth function $f$  on
$\Omega=G(b\Sigma).$

\noindent
(iii) The induced homomorphism $G_*$ of the relative 
homology groups
$$G_* : H_{2n-1} (\Sigma, \Sigma_0;\mathbb R) \mapsto 
H_{2n-1} (G(\Sigma),G(\Sigma_0);\mathbb R)$$
has a nonzero image,
where $\Sigma_0=\Delta \times \partial M.$


Suppose that one of the two cases takes place:

\noindent
(A) $n=1, \ k=2,$

\noindent
(B) $n=2, \ k=3$ and $H_1(M,\partial M)=0.$

Then $Q$ is degenerate on the entire $(2n+1)-$
manifold $\Sigma$, meaning
that $$F=\hat f \circ G$$
on $\Sigma$, for some smooth function
$\hat f$, or, equivalently
$rank\quad dQ = 2n < dim \Sigma=2n+1$
and $ codim _{\ \mathbb R} \ Q(\Sigma)=2.$ The function $f$
is $CR$- function in the interior of $\Omega=G(b \Sigma).$

\end{Theorem}

In the case when $M$ is a closed manifold, i.e.
$\partial M=\emptyset,$
the relative homology groups in condition (iii) of Theorem 2.2
must be replaced by the corresponding
homology groups $H_{2n-1}(\overline\Sigma), \
H_{2n-1}(G(\overline \Sigma)).$

\medskip
Theorem 2.2 can be exposed as a theorem giving conditions of
existence of an extension $\hat f$
of the commutative diagram

\begin{alignat} {6}
&     & b \Sigma  &          &      &                 {}\notag \\
&&{}^G\downarrow & &&  \searrow  {}^F&      &           &   {}\notag\\
&   &\Omega= G(b \Sigma)  &  &&\overset{f}   \longrightarrow &  & \mathbb C
{}\notag.
\end{alignat}
\vskip 0.5cm
to the commutative diagram

\begin{alignat} {6}
&     & \overline \Sigma  &          &      &                 {}\notag \\
&&{}^G\downarrow &  && \searrow  {}^F&      &      &   {}\notag\\
&   &\hat \Omega=G(\overline \Sigma)  &  &&\overset{\hat f}   \longrightarrow &  & \mathbb C  {}\notag,
\end{alignat}
where $F$ and $G$ are $CR-$ mappings.

\bigskip
\noindent
{\bf 2.2. Reduction of Theorems 1 and 2 to Theorem 2.2}

\medskip
Denote $F(\zeta,t)$
the analytic extension of the function
$\zeta \mapsto f(G(\zeta,t))$ from the unit
circle $\{|\zeta|=1\}$ to the unit disc $\{|\zeta| < 1\}.$
Such extension exists
by the condition (*) of Theorems 1 and 2. Then the two functions
$F$ and $G$  define a smooth mapping
$$
Q=(F,G) \ : \overline \Sigma= \overline \Dl \times M \mapsto
\mathbb C \times \mathbb C^n=\mathbb C^{n+1}.$$
\noindent
By definition of the function of $F(\zeta,t)$ the functions
$F$ and $G$ are linked  by the relation

$$
F(\zeta,t)=(f \circ G)(\zeta,t)
$$
\noindent
when  $ (\zeta,t) \in b \Sigma= \partial \Delta \times M.$

\bigskip
\noindent

Now we want to show that Theorem 2.2,  applied to the  functions $(F,G),$
is equivalent to Theorems 1 and 2, applied to the function $f.$

First of all, the conditions (i) of regularity for $G$ are the same
in Theorem 1 and 2, and in Theorem 2.2. 
This condition implies that $G$ is a covering over $\Omega \setminus \partial \Omega$ , i.e. there are finite constant
number of curves through each interior point in $\Omega.$

The condition (ii) in Theorem 2.2 holds because 
$$F(\zeta,t)=f(G(\zeta,t)), \
\zeta \in \partial \Delta, t \in M,$$
by the construction of  $F$. 

Thus, to complete comparing the conditions we need only to check 
that properties (a) in 
Theorems 1 and 2, formulated in terms of analytic discs,
imply the property
(iii) in Theorem 2.2, formulated in terms of the induced homomorphisms of the
homology groups.

\begin{Proposition}
The conditions (a) in Theorems 1 and 2 
is equivalent to condition (iii) in Theorem 2.2.
\end{Proposition}

\begin{proof}

Consider first the case of closed $M$.
Suppose that condition (a) of homological nontriviality holds.
The homology groups of the compact  manifolds
$\Sigma=\overline Dl \times M$ and $M$
are isomorphic and generated by the
the fundamental class $[M]$ which is the homology class of the cycles
$C=\{0\} \times M$ and $M$ correspondingly.

The  $(2n-1)-$ cycle 
$$C =\{0\} \times M$$ 
which we will call {\bf central cycle},
will play an important role in the sequel.

The image cycle $c=G(C)$ intersects each analytic disc
$$D_t= G(\Delta \times \{t\})$$ and by (a) the cycle $c$
represents a nonzero element in $H_{2n-1}(G(\overline \Sigma))$.
Therefore $G_* \neq 0$ and hence the condition (iii) holds.

The case of nonclosed manifolds $M$ and relative homologies  
is considered analogously.

It remains to check that in the planar case, $n=1, \ k=2$ and $\dim M=1,$ 
the condition  of homological nontriviality  converts
to the more transparent condition of Theorem 1 of empty intersection of 
the closures of analytic discs $\overline D_t.$

Let us start with a simple remark.
The domains $\overline D_t$ have a common point,$b,$ 
if and only if the curve $G^{-1}(b)$
intersects each closed disc $$\Delta_t=\overline \Delta \times \{t\}.$$ 
Therefore if the common point $b$
exists then the homology class $[G^{-1}(b)] \in H_1(\overline \Sigma)$ 
is not 0. 

However,  $G$ maps this curve to the point $b,$ and it follows that $G_*=0.$ 
Thus, if (a) in Theorems 1,2 does not hold, then (iii) in Theorem 2.2 fails.   

In the opposite direction, assume that (a) is true, i.e., 
$\cap_{t \in M} \overline{D_t}=\emptyset.$ 
We claim that then the $G$-images of $b\Sigma$ and $\Sigma$ coincide:
$$\Omega=\hat \Omega.$$
Indeed, if not then some value $b$ is taken by $G$ in $\Sigma$ but not
in $b \Sigma.$ Then the function
$$\frac{1}{2\pi i} 
\int\limits_{|\zeta|=1}\frac{dG(\zeta,t)}{G(\zeta,t)-b}$$
is continuous in $t$. It evaluates the number of zeros in the disc $\Delta_t$
and therefore is constant. Since it is different from zero for   least one
value of $t$ it is not zero for all $t$.    Therefore $G$ takes on the value
$b$ on each disc $\Delta$ which means that $b$ belongs to all domains $D_t.$
Contradiction.

Now suppose that (iii) fails. Then $G$ maps the central cycle 
$$C=\{0 \} \times M$$ to a 
cycle $c=G(C),$ (relatively) homologically  
equivalent to zero in $\hat \Omega.$

Consider first the case $\partial M=\emptyset$ and $M=S^1$ in Theorem 1. 
Then the cycle $c$ is homotopic to 0.
Moreover, since the cycle $c$ belongs to the interior of $\hat \Omega=
\Omega,$  it can be contracted, within the interior of $\Omega,$ to a point
$$b \in \Omega \setminus \partial \Omega.$$ 

The mapping $G$ is a   locally trivial foliation
over  $int \ \Omega$ (by the regularity  condition, the
critical values of $G$ are located on
$\partial \Omega=\partial \hat \Omega$).
    Applying the
axiom about covering homotopy (\cite{Hu}, Thm. 4.1.) we  obtain   that
the homotopy $c \sim \{b\}$ lifts up to a homotopy $C \sim C',$ where
$C'$ is a nontrivial 1-cycle  in $\overline \Sigma.$

The new cycle $C'$ projects, by the mapping $G,$ to the point $b$:
$$G(C')=\{b\}.$$
The cycle $C'$ is homotopic to the central cycle $C=\{0\} \times M$
and hence it must intersect each complex disc
$\overline \Dl \times \{t\}$. But this is just another way of saying
that
$$b \in \cap_{t \in M} \overline D_t.$$
Thus, the condition (a) of Theorem 1 fails.

In the case $M=[0,1]$ in Theorem 1, the argument is analogous. 
Namely, if  the condition (iii) in Theorem 1 fails 
and the cycle $c=G(C)$ is relatively trivial, 
then there is a smooth  curve 
$$c_1 \subset G(\Sigma_0)=\overline {D_0} \cup \overline {D_1}$$ such that
$c \cup c_1$ a boundary in the interior of $\hat \Omega= \Omega$. 
The new curve is the image $c \cup c_1 =G(C \cup C_1),$
of the union of the cycle $C$ with  a 1-chain (a curve) 
$$C_1 \subset \Sigma_0= \Delta_0 \times \Delta_1.$$

Then the curve $c \cup c_1$ is homotopic, in the interior of 
$\Omega=\hat \Omega$
to a point. We can perform the homotopy in two steps.

First, we deform $c \cup c_1$ so that the part $c_1$ contracts to a 
point within the set $\overline D_0 \cup \overline D_1.$ 
Then $c$ is deformed homotopically to a closed curve 
(cycle) $c^{\prime} \subset \Omega.$

Again,
by the axiom on covering homotopy, the cycle $C \cup C_1$ 
can be homotopically transformed in $\overline \Sigma$ so that $C_1$ 
contracts to a point within $\Sigma_0$
and the cycle $C$ transforms to a  cycle $C^{\prime}$ 
with end points belonging to $\Sigma_0.$
This new cycle intersects all the discs $\Delta_t.$

Second, we contract $c,$ within the interior of $\Omega,$ 
to an inner point $b \in \Omega.$   
Correspondingly, by the covering homotopy, $C^{\prime}$ can be
transformed to a cycle $C^{\prime\prime} \subset \Sigma,$ 
which intersects each disc $\Delta_t, \ t \in M,$ and
projects by $G$ to the point $b,$ i.e., $$G(C^{\prime\prime})=\{b\}.$$  
This says that for any $t \in M$ 
$$b \in G(\Delta_t)=D_t$$ 
and therefore the condition (a) in Theorem 1 fails. 
The equivalence of conditions (a) and (iii) is proved.

\end{proof}

\noindent

To complete the reduction we have to prove that
the conclusion of Theorem 2.2 implies the conclusions of Theorems 1 and 2.

Indeed, Theorem 2.2 asserts the relation $F=\hat f \circ G$ in $\Sigma.$
This means that $F$ is constant
on the level curves $G=const$ in $\Sigma,$ or, in other words, that the values
of the analytic extensions at a fixed point $ z \in D_t$ do not depend on $t.$
This immediately implies that $f$ is analytic ($n=1$) or $CR-$ function
($n=2$) in $\Omega.$

\medskip
An alternative proof is as follows.
Suppose that $Q$ is degenerate, $F=\hat f \circ G$ on
$\overline \Sigma$ for some smooth function $\hat f$ on $\hat \Omega.$

Let $n=1$. The functions $F$ and $G$ are analytic in $\zeta$ and
differentiating in $\overline \zeta$ the above relation between
$F$ and $G$ yields
$$(\overline \partial \hat f \circ G)
\ \overline{ \partial_{\zeta} G}=0$$
and therefore $\partial \hat f=0$ on $\hat \Omega.$ Thus,
$\hat f$ is analytic on $\hat \Omega$ and hence
$f=\hat f\vert_{ \Omega}$ is analytic on $\Omega.$

Let $n=2$ and the 3-manifold $\Omega$ is locally defined by
the equation $\rho(z_1,z_2)=0,$  where $\rho$ is a smooth real-valued function
in a neighborhood of a point in $\Omega.$ Then
$$\rho(G_1,G_2)=0.$$
Differentiating in the variable $\overline \zeta$
the latter identity and the identity
$$F(\zeta,t)=\hat f(G_1(\zeta,t),G_2(\zeta,t))$$  
leads to a linear system, with the zero determinant due to
the embedding condition $\partial_{\zeta}(G_1,G_2) \neq (0,0).$

Vanishing of the determinant can be written as
$$\overline \partial \hat f \wedge \overline \partial \rho=0,$$
which is just the $CR$-condition for the function $\hat f=f$ in $\Omega.$

This completes the reduction.

\bigskip
\noindent
{\bf 2.3. Necessity of the conditions.}

\medskip
Before going to the proof of Theorem 2.2, let us show necessity of
conditions in Theorems 1,2 and, correspondingly, in their topological versions
Theorem 2.2.

Let us consider the case A, $n=1, \ k=2,$ when the manifold
$\Omega$ is a domain in the plane.

The main condition (a) in Theorems 1 and 2
is that the intersection of the closed analytic discs from the 
family is empty. 
The following example, demonstrating the importance of this condition,
belongs to Globevnik
\cite{G1,G3}. His example is given by function 
$$ f(z)= z^n/\overline z, \ n>1,$$
where $n$ is large enough to provide desirable smoothness at the origin.

The function $\overline z$ extends meromorphicaly
from any circle in the plane inside the disc, with a simple pole at the center,
and with a simple zero, which is inside the disc if and only if
when the point 0 is outside.

It follows that $f$ extends analytically, without poles, inside
any circle enclosing, or even containing, the origin. 
For any family of circles surrounding the origin, the 
corresponding closed discs have the common
point 0, so the condition (a) fails. Nevertheless,
 function $f$ is not holomorphic.

In the context of Theorem 2.2, the above example translates as follows:
$$M=S^1, \quad G(\zeta,t)=R e^{it} + r \zeta,$$
where $R, r$ are nonnegative numbers, 
corresponding to rotation of a circle of radius $r$
 around a circle of radius $R$, and may serve 
a simplest model for understanding the constructions in this article.

The influence of the property of having or nonhaving a fixed point inside all the curves, 
to detecting analiticty in the sense of Theorem 1, 
was  studied in details
by Globevnik in \cite{G3}. He gave there other interesting
examples of the above type. For instance, Globevnik showed that
rotations of equilaterial triangle around its center does not
detect analyticity but this is not true for nonsymmetric triangles.

Define
$$F(\zeta,t)=\zeta e^{it} 
\frac{(R e^{it} + r\zeta)^n}{r e^{it} + R \zeta}, \ n > 1,$$
on the solid torus $\{|\zeta| \leq 1 \} \times S^1.$
In this case the sets $\Omega$ and $\hat \Omega$ are
$$\Omega=G(\partial \Delta \times S^1)=\{|R - r| \leq |z| \leq
R +r  \}, \ \
\hat \Omega=G(\Delta \times S^1)=
\{\max(R - r,0)
 \leq |z| \leq R + r \}.$$

The regularity condition hold as
the Jacobian of $G$ equals to
$$\frac{\partial (G, \overline G)}{\partial (\psi,t)}=
2i R \ r \  sin(t - \psi)$$
and it vanishes if and only if   $\psi-t=0$ or $\pi,$ i.e.
on the circles of radii $R+r$ and $r-R$ centered at 0 which constitute
the boundary of the annulus $\Omega.$

If $r < R$ then $F$ has pole at $\zeta=(-r/R)e^{it} \in \Delta$
and the condition (*) fails.
However, if $R \le r$ then this pole is outside of $\Delta$
and (*) holds. At the same time, the
condition (iii) obviously fails because
the nontrivial cycle 
$$\{(-(R/r) e^{it}, e^{it}) \}$$
is mapped by $G$ to the point $0$ and hence that the induced mapping $G_*$
of the homology groups is trivial. 
Correspondingly, the conclusion of Theorem 1 fails,too,
as $F$ and $G$ are functionally independent in the solid torus and the
mapping $(F,G)$ is nondegenerate.

Consider the case B (Theorem 2), 
related to testing  $CR-$ functions on hypersurfaces.
It is easy to see that the condition (a) can not be removed. 
The corresponding example exists in any dimension. Indeed,
let $\Omega=S^{2n-1}$ be the unit sphere in $\mathbb C^n.$ Consider
the analytic discs which are intersections of the complex unit ball
$B^{2n} \subset \mathbb C^n$ with complex lines $L$ passing through the origin.
The boundaries of the discs are circles,
covering the unit sphere (Hopf foliation).
The condition (a) does not hold, as 0 is the common point of the discs.
Correspondingly, the assertion of Theorem 2 is not true as any function
constant on the the circles $L \cap S^{2n-1}$ but real valued and nonconstant
on the sphere $S^{2n-1}$ provides a counterexample.

\bigskip
\section{Proof of Theorem 2.2}
\noindent
{\bf 3.1. Informal comment on the proof of Theorem 2.2.}

\medskip
The problem under consideration may be viewed as a version of an argument principle for
boundaries of small dimensions (of codimensions greater than 1). Let us explain what do we mean by that.

A simplest one-dimensional analog of Theorem 2.2 could be the following simple fact:
given a smooth mapping $$G:\overline \Delta \mapsto I $$ of the unit
disc onto a segment 
$I=[\alpha,\beta] \subset \mathbb R$, whose restriction to $\partial \Delta$ is a smooth covering over the interval
$(\alpha,\beta)$, and a function
$F$ from the disc-algebra, smooth up to the boundary, then $F=const$ provided $F$
takes equal values on the intersections $G^{-1}(b) \cap \partial \Delta, \ b \in I$.

Here is the proof in the form we  need (a similar argument was used in \cite{BW}, \cite{AV}, \cite{G1}).
Suppose there exist $b \in F(\overline \Delta) \setminus F(\partial \Delta).$
Then by the argument principle
$$\frac{1}{2\pi i} \int\limits_{\partial \Delta}\frac{dF}{F-b} = \#
\{z \in \Delta: F(z)=b  \} > 0.$$
On the other hand, from the condition we have
$F=f \circ G$ on $\partial \Delta$ for some smooth function $f$ and
then
$$\frac{1}{2 \pi i} \int_{\partial \Delta} \frac{dF}{F-b}=\deg G \cdot
\frac{1}{2 \pi i}
\  \int\limits_{I}\frac{df}{f-b}=0,$$
because the Brouwer degree of the mapping 
$G:\partial \Delta=S^1 \mapsto
I$ is zero, $\deg G=0.$

The obtained contradiction says that
$F(\overline \Delta)=
 F(\partial \Delta)$ which implies $f=const$, as
otherwise, the curve $F(\partial \Delta)$
must contain an open set, due to the openness property.

The above argument can be viewed as a simple model
of the idea of the proof of our main result, Theorem 2.2. However, the main
difference and the main difficulty is that in our case the function $F$ depends
on the additional parameter $t$ and the fibers of the foliation $G$ do not belong
to the complex discs where $F$ is analytic. Moreover, the fibers are transversal
to the discs.
Our goal is to prove that then, still, $F$ is constant on the $G$-fibers and
hence the mapping $Q=(F,G)$ is degenerate: $Q$ maps the level curves $G^{-1}(b)$ to points. 

The topological interpretation of the above toy model might be the following:
if the analytic function $F$ collapses the boundary circle to a curve
then $F$ collapses the unit discs as well, and the images of the disc
and its boundary coincide:
$F(\overline \Delta)=F(\partial \Delta).$

Theorems 2.2 may be interpreted in the similar topological terms.
Take  for example, the case
$n=1, \ k=2$, corresponding to the strip-problem for a periodic family
of Jordan curves.

The mapping $G$ from Theorem 2.2 maps the 2-dimensional torus $T^2$ 
to a planar domain, $\Omega.$ Now, the function
$F$ is constant on the level sets of $G$, as $F=f \circ G$, and therefore the
composite mapping $Q$, defined on the solid torus,
$$Q=(F,G): \Sigma \mapsto \mathbb C^2 \cong \mathbb R^4$$ collapses the boundary torus $T^2$
to a manifold diffeomophic to a planar domain (which is the graph of the function $f$ over $\Omega.$)

The image $Q(\Sigma)$ of the interior is a set in $\mathbb R^4$ attached to the 2-dimensional
image of the boundary torus $T^2$, and this image is ``flat'', i.e. is tolopologically equivalent to
a planar domain. 
In Theorem 2.2 we prove in that if the mapping $Q$ is homologically nontrivial, meaning that 
$Q$ induces nontrivial homomorphism of the homology group, then 
$Q$ collapses the interior of the solid torus as well it does with its boundary. In other words, the image of the 3-dimensional
solid torus $\Sigma$ is 2-dimensional  and is 
contained in the image of the boundary and $Q(\overline \Sigma)=Q(\partial \Sigma),$ similarily to our model example.

\bigskip
\noindent
{\bf 3.2. The plan of the proof of Theorem 2.2.}

\medskip
We assume, to the contary of the assertion, that $Q$ is nondegenerate 
and our goal will be arriving to a contradiction.

\bigskip
\noindent

{\it Step 1.}

\medskip
A key point is Lemma 3.3 where we prove the  symmetry relation (2) of
linking numbers which are
periods (winding numbers) of one function 
on the zero sets of the other.
Lemma claims that certain logarithmic residue type integrals are equal. This equality is true for functions analytic in $\zeta$ variable
and having the important property: they take  equal values at the points on the boundary manifold $b \Sigma$, belonging to the same
$G$-fiber.

Lemma 3.3 is proved by the standard technique from residues theory,
assuming removing neighborhoods of singularities, applying Stokes formula to
differential forms on the remaining manifold and then shrinking the removed
neighborhoods.

The key point is that the constancy of functions on the $G$- fibers on $b \Sigma$
provides cancellation of a surface term in the Stokes formula and this
leads to the needed symmetry relation.

\bigskip
\noindent
{\it Step 2.}

\medskip
In the second, analytic, part of the proof we study the
minors $J$ of the  Jacobi matrix $J(F,G)$ in $\Sigma$.
First, we observe that the zero sets of $J$ (the $Q$-critical sets)
contain the central cycle $C=\{0\} \times M.$  Since we assumed that
$Q$ is not degenerate, there are minors $J$ which 
are not identically zero and hence 
the zero sets $J^{-1}(0)$ are essentially  $(2n-1)$-dimensional.

The final part of the proof, for instance in
the case $\partial M=\emptyset,$  (Theorems 1 and 3), is as follows.

Take for simplicity the case $n=1.$
We apply the symmetry relation proved in Lemma 3.3 to the functions $J$ and $G-b$. However, Lemma
3.3 requires that both functions identify points on $b \Sigma$ from the same $G$-fiber. The function $G-b$
certainly satisfies this property, but $J$ may not. Nevertheless, it appears that the ``phase'' part, $\Theta=J/\overline J$,
does satisfy the required property away from zeros of $J$. This is enough to write the symmetry relation
for $J$ and $G-b$, as the change of argument of $\Theta$ is just the doubled change of argument of $J.$

Now we use the symmetry relation to obtain contradiction with the assumption that $J$ does not vanish identically.

First, take $b \notin \hat \Omega=G(\overline \Sigma)$. Then $G^{-1}(b)=\emptyset$
and the integral over $G^{-1}(b)$ in the symmetry relation (2) vanishes because the set of integration is empty. 

On the other hand, condition (iii) guarantees that
the image  $c=G(C)$ of the nontrivial loop $C=\{0\} \times M$, is homologically nontrivial in $\hat \Omega$. 
Since $C \subset J^{-1}(0)$ then  $c \subset G(J^{-1}(0))$.

Then $b$ can be choosen out of $\hat \Omega$ and so that the loop $c$ has a nonzero index, $ind_b c \neq 0$, 
with respect to the point $b$
and we have the desired contradiction because then      
the second integral in (2), over $G(J^{-1}(0))$, is different from zero.

This contradiction implies that $J=0$ identically
which is just the         conclusion of Theorem 2.2.

In the case  $\partial M \neq \emptyset$ (Theorems 2 and 4) we slightly modify
the above argument. First, an additional surface term
appears in Stokes formula, corresponding the the boundary of $M.$
Then the symmetry relation delivers
a jump-function which counts algebraic number of intersections
with the image $G(J^{-1}(0))$ of the  critical set.

The homological condition (iii) 
provides that there exists a path, $L,$
which crosses the union of the discs $D_t, t \in M$ without meeting
the discs with $t \in \partial M$ and has the nonzero 
intersection index with $G(C) \subset G(J^{-1}(0))$.   

The final argument rests on determining the relative
homological class of a cycle of codimension 1 by intersection indices
with a transversal curve $L$.
The contradiction comes from
comparing the total algebraic number of the intersection indices
(jumps of the counting function) which is different from 0 with the total variation of
the counting function along the path $L$ which appears to be 0.

\bigskip
\noindent
{\bf 3.3. The Brouwer degree of the mapping $G.$}

\medskip
Let us study properties of the mapping $G.$

We start with
the planar case $(2n,k)=(2,2)$ and with the case of the closed manifold $M$. In this case
$$M=S^1, \ \Sigma=\Delta \times S^1, b \Sigma = \partial \Sigma
=S^1 \times S^1=T^2.$$
\begin{Lemma}
The Brouwer degree of the mapping
$$G : T^2 \mapsto G(T^2)=\Omega $$
is zero, $\deg G=0.$
\end{Lemma}

\begin{proof} The image $G(T^2)$ is a compact domain
 $\Omega=G(T^2) \subset \mathbb C.$

Embed $\Omega$ to the Riemann sphere as a compact subset
$\tilde \Omega \subset S^2.$ Then we can view the mapping $G$ as a mapping to
$S^2$:
$$ G: T^2 \mapsto \Omega \subset \tilde \Omega \subset S^2.$$

The 2-dimensional homology groups of the 2-dimensional
compact manifolds without boundaries are
$$ H_2(T^2; \mathbb Z) = H_2(S^2; \mathbb Z)= \mathbb Z$$
and the Brouwer degree is defined by the relation
$$ G_* (\mu_1)=\deg G \ \mu_2,$$
where $\mu_1, \mu_2$ are generators (fundamental classes) of the
corresponding homology groups 
$H_2(T^2; \mathbb Z)$ and $H_2(S^2;\mathbb Z).$  
However, the image of $G$ is a compact
subdomain of the sphere $S^{2}$ and is contractible, hence
$G_*(\mu_1)=0$ and therefore $\deg G=0.$
\end{proof}

We will need a local version of Lemma 3.1:
\begin{Corollary}
If $O \subset \Omega$ is a submanifold, then the local Brouwer
degree of the restricted mapping
$$G: T^2 \setminus G^{-1}(O) \mapsto \Omega \setminus O$$
is zero.
\end{Corollary}
\begin{proof} Recall that locally the Brouwer degree is defined
as the algebraic  number of points in the
preimage $G^{-1}(b)$ of a regular
value $b \in \Omega,$ i.e.
$$ \deg G=\sum_{G(p)=b} sign \  \ dG_p.$$
Clearly, this number preserves after simultaneous removing a set
$O$ from  $\Omega =G(T^2)$ and its
full $G-$ preimage $G^{-1}(O)$ from $T^2.$
\end{proof}

\bigskip
\noindent
{\bf 3.4. Linking numbers.}

\medskip
Everywhere in this section we assume that $\partial M=\emptyset,$ but
combine both cases $(2n,k)=(2,2)$ and $(2n,k)=(4,3),$
so that $\Omega$ is either a planar domain or
3-dimensional closed manifold in $\mathbb C^2.$
Note that if $\partial M=\emptyset$ then
$b \Sigma$ coincides with the topological boundary,
$b \Sigma =\partial \Sigma.$

Denote $\beta_{MB}$ the Martinelli-Bochner $(2n-1)-$
differential form $\beta_{MB}$
(Martinelli-Bochner kernel) in the space $\mathbb C^n$ (see,e.g.
\cite{AYu},\cite{St1}).
It will be convenient to denote the coordinates in $\mathbb C^n$ by
$z_2, \cdots, z_{n+1}.$ Then
$$ \beta_{MB} = \frac{1}{(2\pi i)^n} |z|^{-2n}\omega^{\prime} (\overline z)
\wedge \omega (z),$$
where $\omega^{\prime}(\overline z)=
\sum_{j=2}^{n+1}(-1)^j \overline z_j \ d\overline z[j]$ and
$\omega(z)= dz_2 \wedge\cdots\wedge dz_{n+1}.$
and $[j]$ means that the differential $d\overline z_j$ is skipped.

The Martinelli-Bochner form coincides with the $(2n-1)-$
surface form in the real Euclidean space $R^{2n}$ and integration in
this form is well related with computing degrees of mappings, and
linking numbers (see e.g., in the close context,\cite{AW},\cite{St1}.)

Define the $2n$-form in $\mathbb C^{n+1}$:
$$\beta = \frac{1}{2\pi i} \frac{dz_1}{z_1} \wedge \beta_{MB}.$$
In particular, for $n=1$ we have
$$\beta=\frac{1}{(2 \pi i)^2} \frac{dz_1}{z_1} \wedge \frac{dz_2}{z_2}.$$

The form $\beta$ is closed in
$\mathbb C^{n+1}$ with the deleted linear spaces $\{z_1=0\}$ and
$ \{z^{\prime}=(z_2, \cdots,z_{n+1})=0\}.$

The following lemma is a key one in the proof of Theorem 2.2.
It asserts a symmetry of
linking numbers associated to
the level sets of functions, compatible with the foliation
$G.$ By {\it $G-$ compatibility}
we understand the constancy on the $G-$ fibers.

\begin{Lemma} Let $G:\overline \Sigma \mapsto \mathbb C^n$
be the mapping from Theorem 2.2, and $\partial M=\emptyset.$
Let $J$ be a nonzero $C^1-$
function
$$ J: \overline \Sigma \mapsto \mathbb C$$
having the following properties:

\noindent
1. The zero set $J^{-1}(0)$ is a $(2n-1)-$ chain in
$\overline \Sigma.$

\noindent
2. The set $J^{-1}(0) \cap b\Sigma$ contains the critical set $Crit(G)$
of the mapping $G$ on $b\Sigma.$

\noindent
3. The quotient $$\Theta =\frac{J}{\overline J}= exp(2i \ arg \ J):
(\partial \Delta \times M) \setminus J^{-1}(0) \mapsto S^1$$
is $G$-compatible, i.e., $\Theta$
is representable, outside of zeros of $J,$ in the form

\begin{equation}
\Theta = \sigma \circ G,
\end{equation}
for some smooth function $\sigma$ defined on $G(b\Sigma \setminus J^{-1}(0)).$

Then

\noindent
a) for any regular value $b \in \mathbb C^n$ of the mapping $G,$
such that $J^{-1}(0) \cap G^{-1}(b) = \emptyset,$ holds:
\begin{equation}
2 \int_{J^{-1}(0)} (G-b)^*\beta_{MB} =
\frac{1}{2\pi i} \int_{G^{-1}(b)}
(\frac{dJ}{J}-\frac{d\overline J}{\overline J}).
\end{equation}

\noindent
b) The chain $c=G(J^{-1}(0))$ is a cycle, i.e. $\partial c =\emptyset.$

\end{Lemma}

\begin{proof}
We use the usual technique from the theory of residues, which involves
deleting neighborhoods of singular sets, applying Stokes formula to the
remaining manifold and then
shrinking the neighborhoods (e.g. \cite{AW}, \cite{St1}).

\bigskip

First of all, we will complete the set $J^{-1}(0)$ to 
make it $G$-compatible on $b \Sigma,$ i.e. union of $G$-fibers.
For this aim, consider the full $G$-preimage 
$$N_0=[G^{-1}(G(J^{-1}(0)))] \cap b\Sigma. $$  Define $N$ the compliment of the zeros of $J$ on $b\Sigma$ in $N_0:$ 
$$N= N_0 \setminus (J^{-1}(0) \cap b\Sigma).$$
By the construction, the  union $$(J^{-1}(0) \cap b\Sigma) \cup N=N_0$$ 
is $G$-compatible on $b \Sigma.$ 
The set $N$ is either empty,
or, by the regularity conditions for $G,$ 
it is finite union of smooth manifolds of the dimensions
less or equal to $2n-1.$  

The first step of the proof of the formula (2) is to 
construct a family of shrinking
neighborhoods, in $\Sigma,$ of the singular sets $J^{-1}(0)$
and $G^{-1}(b)$. 

Define  $$A_{\varepsilon}=\{|J| < \varepsilon \}, \
B_{\varepsilon}= \{|G-b| < \varepsilon \}.$$ 
Additionally, construct
a family $N_{\varepsilon}$ of $\varepsilon$-
neighborhoods of the set $N.$

Now remove the constructed neighborhoods from 
$\Sigma$ and denote the remainder:
$$\Sigma_{\varepsilon}= \Sigma \setminus
(A_{\varepsilon} \cup B_{\varepsilon} \cup N_{\varepsilon}).$$
Define the orientation on
 the manifolds $\partial A_{\varepsilon}, \ \partial B_{\varepsilon}$
and $\partial N_{\varepsilon} $
by the inward, with respect to $\Sigma_{\varepsilon},$ normal vector.

Define on $\overline \Sigma_{\varepsilon}$ the $2n$-form:
$$\Xi = \Phi^* \beta = \frac{1}{2 \pi i}
\frac{d\Theta}{\Theta} \wedge (G-b)^*\beta_{MB},$$
where $\Phi=(\Theta,G-b).$
The differential form $\Xi$ has no singularities in $\overline
\Sigma_{\varepsilon}$
and is closed there because the form $\beta$ is closed.

Therefore the Stokes formula yields:
$$
\int_{\partial \Sigma_{\varepsilon}} \Xi =
\int_{\Sigma_{\varepsilon}} d\Xi = 0.
$$
The boundary of $\Sigma_{\varepsilon}$ consists of the four parts equipped
by the induced orientations:
$$\partial \Sigma_{\varepsilon}=\partial A_{\varepsilon}
\cup \partial B_{\varepsilon} \cup \partial N_{\varepsilon} 
\cup (b \Sigma \setminus
(A_{\varepsilon} \cup B_{\varepsilon} \cup N_{\varepsilon}))$$
and therefore the latter identity reads as
\begin{equation}
\int_{b \Sigma \setminus (A_{\varepsilon} \cup B_{\varepsilon} 
\cup N_{\varepsilon})} \Xi
-\int_{\partial A_{\varepsilon}} \Xi -
\int_{\partial B_{\varepsilon}} \Xi -\int_{\partial N_{\varepsilon}}=0.
\end{equation}
\bigskip
\noindent

The signs minus before the  last three integrals appear because
the induced orientations on $\partial A_{\varepsilon}$ and
$\partial B_{\varepsilon},$ required by Stokes formula, are opposite
to the above defined orientations.

Now we  take the limit
of the integrals, when $\varepsilon \to 0.$

First, let us check that the added set $N$ contributes nothing, i.e. the integral over $\partial N_{\varepsilon}$
in (3) dissappears when $\varepsilon \to 0$
Indeed, since $N$ is not contained in $J^{-1}(0)$ 
then by real-analyticity
$\dim (\overline N \cap J^{-1}(0)) < 2n-1 .$  
On the other hand, the intersection  $\overline N \cap G^{-1}(b)$ with the level curve of $G$ is not more than finite. 

The differential form  $\Xi$ is the wedge product of the Martinelli-Bochner type form and the logarithimic derivative type form.
The first factor has singularities at zeros of functions $J$ and the second one-on the curve $G=b.$ 
We see from the dimensions of these singularities 
on $N$ that $\Xi$ is integrable (in a regularized sense) on $N$ 
and since $\dim N < 2n, \ \deg \Xi=2n,$ 
the integral of $\Xi$ over $N$ is zero.
Since  the neighborhoods  $N_{\varepsilon}$ shrink to $N,$ 
the last integral in the left hand side in (3) goes to zero as
$\varepsilon$ goes to zero.

To understand the limits of other integrals in (3), write
the differential form under the integral (3) as
$$ \Xi= \frac{1}{2 \pi i} (\frac{dJ}{J}-\frac{d\overline J}{\overline J})
\ \wedge (G-b)^*\beta_{MB}=
\Phi^* \beta - \overline \Phi^*\beta,$$ where $\Phi=(J, G-b),
\overline \Phi=(\overline J, G-b)$ and
$$\beta=\frac{1}{2\pi i}\frac{dz_1}{z_1} \wedge \beta_{MB}.$$

The mapping $\Phi$ maps
$J^{-1}(0)$ to $\{0\} \times (G-b)(J^{-1}(0))$ and hence
change of variables  $z=\Phi (\zeta,t) \in \mathbb C^{2},$ i.e.
$z_1=J, \ z_2=G-b,$ in (3) leads, after
letting $\varepsilon \to 0,$ to:

\begin{equation} 2 \int_{J^{-1}(0)} (G-b)^*\beta_{MB}-
\frac{1}{2\pi i} \int_{G^{-1}(b)} (\frac{dJ}{J}-
\frac{d\overline J}{\overline J}) = \int_{b \Sigma^{\prime}} \Xi,
\end{equation}
where 
$$b\Sigma^{\prime}= 
b\Sigma \setminus (J^{-1}(0) \cup N \cup G^{-1}(b))=
b\Sigma \setminus (N_0 \cup G^{-1}(b)) .$$

Here we  have used that  the
forms $(2 \pi i )^{-1} dz_1/z_1$
and  $\beta_{MB}$ are correspondingly the arc and the surface measures
on the unit circle and the unit sphere
in the spaces of the variables $z_1$ and $z^{\prime}.$ 
The sign minus before the second integral
comes from the orientation consideration. The factor 2 in the left hand side
is due to the relation
$$\frac{d \overline J}{\overline J}=-\frac{dJ}{J}$$
on $|J|=\varepsilon.$

Formula (2) will be proved if we would show that the right hand side in (4)
is zero.
We will prove this separately for the case when $\Omega$ is a 
domain in $\mathbb C$ (the case A, (2n,k)=(2,2)) and
for the case when $\Omega$ is a hypersurface in $\mathbb C^2$ 
(the case B, (2n,k)=(4,3)).

In the first case, 
 $(2n,k)=(2,2)$ we use Corollary 3.2 from
Lemma 3.1. By regularity condition , the mapping $G$
is a covering of $b \Sigma \setminus Crit(G)$
over $\Omega \setminus \partial \Omega$.

The set $b \Sigma ^{\prime}$ 
is $G$-compatible because it is obtained by removing from $b \Sigma$
the full $G$-preimages,$N_0$ and $G^{-1}(b).$ It does not contain the 
critical $Crit(G)$ because this set is removed along with $J^{-1}(0).$
Hence $G$ is a covering on $b \Sigma^{\prime}.$
By Corollary 3.2 from Lemma 3.1, we have
that $$\nu=\deg G\vert_{b \Sigma^{\prime}} =0.$$

By the condition 2 ,
$\Theta=\sigma \circ G$ away from  zeros of $J.$ Then the change of
variables $G(\zeta,t)=z^{\prime}$ in the integral yields:

$$\int_{b \Sigma^{\prime}}\Xi=
\nu \cdot \frac{1}{2\pi i}
\int_{G(b\Sigma^{\prime})} 
\frac{d\sigma}{\sigma} \wedge \beta_{MB}(z^{\prime}-b) =0,$$
because $\nu=0.$

For $(2n,k)=(4,3),$ the vanishing of the integral on the right hand 
side of (3) follows
even more simply. In this case, the mapping $G: b \Sigma \mapsto \Omega$
is a smooth foliation of 4-dimensional manifold over a
3-dimensional manifold, with one-dimensional fibers. 
The differential form 
$$\eta=\frac{d\sigma}{\sigma} \wedge \beta_{MB}(z^{\prime}-b)$$
has degree 4 and therefore $\eta=0$ on the 3-dimensional manifold $\Omega.$
Then the pull back of the form $\eta$ to $b \Sigma$ satisfies
$\Xi=G^*(\eta)=0$ and therefore the integral of $\Xi$ over 
$b\Sigma^{\prime} \subset b\Sigma$ vanishes.

\bigskip
Now we have completed the proof of the main part, a), of Lemma 3.3. 
Notice that 
the above proof may be briefly exposed by using the language
of currents (see, \cite{AW}).

It remains to check  the assertion b). To this end,
take in (2) the point $b$ belonging to the unbounded component
$V$ of $\mathbb C^n \setminus G(\overline \Sigma)$. Then the right hand side
in (2) equals zero and hence, after changing
of variables $z^{\prime}=G(\zeta,t)$ in the left hand side, we obtain
$$\int_c \beta_{MB}(z^{\prime} -b)=0, \ b \in V.$$
It can be proved by many ways that vanishing of Martinelli-Bochner type integral in a neighborhood 
of infinity implies that the surface of integration is closed, $\partial c=0.$

Now the proof of Lemma 3.3 is completed.

\noindent
{\bf Remarks.} 

\noindent
1.The proof of Lemma 3.3  may be 
briefly exposed by using the language and the technology
of currents (see, \cite{AW}).

\noindent
2. In the proof of Theorem 2.2 for case of closed manifold $M$ we will
use formula (2) in the situation when $G^{-1}(b)=\emptyset$ and therefore
the right hand side in (2) is zero. 

\end{proof}

\bigskip
\noindent
{\bf 3.5. The Jacobi determinants $J$, the function $\Theta,$
and the central cycle $C,$ for the case (2n,k)=(2,2).}

\medskip
Our immediate goal is to construct a functions $J$ and
satisfying the conditions
of Lemma 3.3, and vanishing on the central cycle $C=\{0\} \times M.$
The function we are going to construct will also carry information about degeneracy of the mapping $Q.$

Such function turns out to be the Jacobian of the mapping
$Q=(F,G)$ in the angular variable $\psi$ and the local coordinate $t$
on $M.$
The needed properties follow from the functional relation between
$F$ and $G$ on $b \Sigma,$ while
vanishing  on the central cycle $C= \{\zeta=0\}$ in
$\overline \Sigma $ comes from the vanishing
at $\zeta=0$ of the tangent vector field to the unit circle. This vector field acts on
analytic functions as the complex derivation
$\partial_\psi= i\zeta \partial_{\zeta}.$

Let us do the corresponding
computations.
Choose the basis $\partial_\psi, \partial_{t}$
in the tangent space to $b \Sigma.$
Here $\zeta=e^{i\psi}$ and $t$ is the local coordinate on $M.$
If $M=[0,1]$ then take $t \in [0,1]$  and if $M=S^1$
then $t$ can be taken the angular variable on the unit circle,
$M=S^1=\{e^{i t}\}.$

We will use the notation $\nabla$ for the column:
$$\nabla=column \ [\partial_\psi,\partial_{t}].$$
In this notation, the Jacobi matrix for the mapping $G$
becomes $[\nabla G, \nabla \overline G].$

\begin{Lemma} Define
$$J_+ =
\det \ [\nabla F_, \nabla \overline G],$$
$$J_{-} =
\det \ [\nabla G, \  \nabla F].$$

Thus, $J_{+}$ is obtained from the Jacobi matrix of $G$ by replacing
$\nabla G$ by $\nabla F$, while $J_{-}$ - by  replacing
$\nabla \overline G$ by $\nabla F.$

Then on $b \Sigma$ the relation holds:

\begin{equation} J_{+} = \det \  [\nabla G, \nabla \overline G] \
(\partial f \circ G), \ \ \
J_{-}=\det \  [\nabla G, \nabla \overline G] \
(\overline \partial f \circ G),
\end{equation}

Here $\partial, \overline \partial $ are the derivatives in
$z$ and $\overline z$ respectively.
\end{Lemma}

\begin{proof}

We start with the link $F=f \circ G$
between $F$ and $G$ (condition (ii) of Theorem 2.2)
on the manifold $b \Sigma=\partial \Delta \times M$.
Differentiation in local coordinates on $b \Sigma$
and the chain rule lead to the linear system:
\begin{equation}
[\nabla G, \nabla \overline G]  \ \
[\partial f \circ G, \ \overline \partial f \circ G] = \nabla F.
\end{equation}

Then (5) follows immediately from the Cramer's rule when solve the
linear system (6) for $[\partial f \circ G, \ \overline \partial f \circ G].$
\end{proof}

In the sequel we will be exploiting only the ``minus'' minor, $J_{-},$
as it possesses needed orientation properties.
The function $J(\zeta,t)=J_{-}(\zeta,t)$ is just the Jacobian
$$J = \frac {\partial (F,G)}{\partial (\psi, t)}=
\partial_{\psi}F \ \partial_t G- \partial_{\psi}G \ \partial_t F,$$
where $\zeta=|\zeta|e^{i\psi}.$

On the manifold $S^1 \times M,$ the function $J$
also can be understood as the Poisson bracket $J=\{F,G\}.$
The function $J$ is defined in the entire $\Delta \times M$ and can be
expressed there in terms of complex
derivatives in $\zeta:$
\begin{equation} J(\zeta,t)=i\zeta ( \partial_{\zeta}F \partial_t G -
\partial_{\zeta} G \partial_t F).
\end{equation}
Consider a smooth $G-$ level curve $\Gamma=G^{-1}(b),$
parametrized by  $\zeta=\zeta(t).$
Differentiating  the identity $G(\zeta(t),t)=b$
in $t$ and taking into account that $\partial_{\zeta}G \neq 0$
we obtain for the directional derivative $\partial^G$ along the 
$G$- level curves:
$$\partial^{G} F=\frac{d}{dt}F(\zeta(t),t)=
\partial_{\zeta} F  \ \zeta^{\prime} (t) + \partial_t F=
- \frac{\partial_{\zeta} F \ \partial_t G -
\partial_{\zeta}G \ \partial_t F}{\partial_{\zeta} G}$$
and therefore $J$ is related to the directional derivative by
$$J=-i \zeta \partial_{\zeta}G \ \partial^G F.$$

Thus, vanishing of $J$ identically in $\Sigma$ means that $F=const$
on the $G-$ level curves and therefore $F$ is a function of $G.$

\begin{Lemma} The functions $J=J_{-}$
possesses the properties:

\noindent
1. $J(0,t)=0$ for all $t \in M.$

\noindent
2. The zero set $J^{-1}(0)$ contains the critical
set $Crit(G)= \{ det[\nabla G, \nabla \overline G]=0 \} 
\subset b\Sigma.$

\noindent
3. On $b \Sigma,$ away from zeros of $J$,  the relation holds
$$\Theta=\frac{J}{\overline J}=\sigma \circ G,$$
where
$$\sigma=-\frac{\overline \partial f}{\overline  {\overline{\partial} f}}.$$

\end{Lemma}

\begin{proof}

The property 1 follows from (7). 
The property 2 follows from the representation of $J$ in formula (5). 
The property 3 also follows from (5)
because the Jacobian $\det [\nabla G,\nabla \overline G]$
obviously takes purely imaginary values.
\end{proof}

Thus, the function $J$ satisfies all the conditions of Lemma 3.3, except, 
maybe, the assumption 1 about  a nice structure of the zero set $J^{-1}(0).$

\bigskip
\noindent
{\bf 3.6. The functions $J$ in the case (2n,k)=(4,3).}

\medskip
\noindent
{\it Constructing  Jacobi minors.}

Now we need a version of Lemma 3.5 for the case B, when $n=2$ and 
$\Omega$ is
a hypersurface in $\mathbb C^2, \ \ \dim \Omega=2n-1=3.$
Recall that in this case, the regularity condition says that
$rank \ dG\vert_{\partial \Delta \times M} = 3.$

Consider the complex 
Jacobi matrix $$[\nabla G_1,\nabla G_2,\nabla F]$$
contains 3 columns and 4 rows. The rows correspond to the derivatives
with respect to $\psi, t_1,t_2,t_3,$ where $t_k$ are local coordinates
on $M.$

Let us first specify how should we understand the differenetiation in 
local coordinates  $t_k$.

Any 3-dimensional compact orientable manifold $M$ is 
parallelizible (the tangent bundle is trivial), 
see e.g.\cite{MS}. This means that 
there exist three linearly independent smooth
tangent vector fields $X_1,X_2,X_3 \in T(M)$ on $M$. 

In a neighborhood of any point in $M$ we define parametrization $t_k$
via coordinates in the basic $X_1,X_2,X_3$ in the tangent plane, so that
the vector field $X_k$ corresponds to the differentiation in $t_k$,
$X_k=\partial_{t_k}, k=1,2,3.$

The degeneracy of the mapping $Q=(F,G)$ would be shown
if we prove that  all $3 \times 3$ minors
$J^0, J^1,J^2,J^3 $ are identically zero ,
where $J^0$ is obtained from the above Jacobi matrix
by removing the first row (containing the derivaties with respect to
$\psi$)
and $J^k,k=1,2,3,$ are obtained by deleting the row
$(\partial_{t_k}G_1,\partial_{t_k}G_2,\partial_{t_k}F).$ 
\begin{Lemma} $J^1,J^2,J^3 \equiv 0$ imply $J^0 \equiv 0.$ 
\end{Lemma}
\begin{proof} 
Denote
$$h_{\psi}=\partial_{\psi}(G_1,G_2,F), 
h_{\alpha}=\partial_{t_{\alpha}}(G_1,G_2,F),\alpha=1,2,3.$$
Suppose $J^0(\zeta_0,t_0) \neq 0.$ Then $J^0(\zeta,t) \neq 0$
for $(\zeta,t)$ in a neighborhood of $(\zeta_0,t_0)$ and 
the vectors $h_1,h_2,h_3$ are linearly independent over $\mathbb C$
for such $(\zeta,t)$.

On the other hand, the condition implies that each system consisting of
the vector $h_{\psi}$ and 
any two vectors from the set  $h_1,h_2,h_3$, is 
$\mathbb C$-linearly dependent.
Since the last two vectors 
in the above triple are linearly independent, we have
$$h_{\psi} \in \Pi_{\alpha,\beta}=span\{h_{\alpha},h_{\beta}\},
\alpha,\beta=1,2,3, \alpha \neq \beta.$$ 

But due to the linear independence of the vectors
$h_1,h_2,h_3$ we have
$\Pi_{1,2} \cap \Pi_{1,3} \cap \Pi_{2,3} = \{0\}.$
The vector $h_{\psi}$
belongs to this intersection and hence
$$h_{\psi}=(\partial_{\psi}G_1,\partial_{\psi}G_2,\partial_{\psi} F)=
(0,0,0).$$ 
The derivatives are evaluated at points $(\zeta,t)$
in an open set in $\Sigma.$. However it is impossible as the functions 
$G_1(\zeta,t),G_2(\zeta,t)$
are nonconstant analytic (in $\zeta=re^{\psi}$) functions.
\end{proof}

Lemma shows that, to prove the degeneracy of the mapping $(F,G)$,
it  suffices to check vanishing  the minors $J^1,J^2,J^3$.
The advantage of these minors is that their first
row is $\partial_{\psi}(G,F)$ and therefore vanishes at $\zeta=0,$ the fact
which is crucially important to us. The minor $J^0$ does not have this
property.

Let us fix some of the minors, for instance, $J=J^3.$ This minor
corresponds to  deleting from Jacobi matrix the last row with the derivatives  
with respect to  $t_3:$
$$J=
\vmatrix \partial_{\psi} G_1 & \partial_{\psi}G_2  & \partial_{\psi} F \cr
\partial_{t_1} G_1 & \partial_{t_1} G_2 & \partial_{t_1} F \cr 
\partial_{t_2} G_1 & \partial_{t_2} G_2 & \partial_{t_2} F 
\endvmatrix 
$$
Remind that the local parameters on $M$  are chosen so that $X_1$ and $X_2$
are differentiations in the first two coordinates
$X_1=\partial_{t_1},X_2=\partial_{t_2}.$

To the end of this section,in order not to overload notations, 
we will denote by the symbol ``$\det$''
the minors including the rows with the derivatives
$\partial_{\psi},\partial_{t_1},\partial_{t_2}$, but not $\partial_{t_3}.$

\bigskip
\noindent
{\it Properties of the function $J.$}

We have chosen one of the minors $J^1,J^2,J^3$, namely $J=J^3$,
but what follows applies for any minor $J^{\alpha}.$

Suppose that $J=J^3$ is not identically zero.
We start with  checking that $J/\overline J$ is $G$-compatible on $b\Sigma$
away from the zeros of $J.$

It would be more convenient to perform computations for arbitrary $n$
and then set $n=2.$
Let $\rho=\rho(z,\overline z)$
be a smooth real valued 
defining function in a neighborhood $U$ of the hypersurface
$ \Omega:$
$$\Omega=\{ z \in \mathbb C^n :\rho(z, \overline z)=0\}, \
\nabla \rho \neq 0.$$
The function $f$ is assumed to be extended as a smooth function in $U.$

Consider the tangential $CR-$ operators on $\Omega:$
$$\overline \partial_{\mu,\nu}= \overline \partial_\mu \rho \ \overline
\partial_\nu -
\overline\partial_\nu \rho \ \overline \partial_\mu.$$
Here $\partial_\mu$ and $\overline \partial_\nu$ are derivatives in $z_\mu$
and $\overline z_\mu$ correspondingly. The system 
$\overline \partial_{\mu,\nu}, \ \mu < \nu,$
forms a basis in the complex tangent $\overline \partial$-bundle 
$T^{0,1}(\Omega).$

We start again with the main relation on the boundary
manifold $b \Sigma= \partial \Delta \times M:$
$$ F=f \circ G $$
and differentiate it in the local coordinates 
$\psi,t:$
\begin{equation}
\sum_{j=1}^n (\partial_j f \circ G) \nabla G_j +
\sum_{j=1}^n (\overline\partial_j f \circ G) \nabla \overline G_j=\nabla F.
\end{equation}
The extra relation between the gradients
comes from the equation of the hypersurface $\Omega:$
$$\rho(G_1,\cdots,G_n,\overline G_1,\cdots,\overline G_n)=0.$$
Differentiating yields:
\begin{equation}\sum_{j=1}^n (\partial_j \rho \circ G) \nabla G_j+
\sum_{j=1}^n (\overline {\partial}_j \rho \circ G) \nabla \overline G_j=0.
\end{equation}
Fix $\mu$ such that  $\overline \partial_\mu \rho \neq 0$ in a neighborhood
of a point in $\Omega.$
Express
$\nabla \overline G_\mu$ from (9) and substitute in (8). After grouping terms
we obtain:
$$\sum_{j=1}^n [(\partial_{\mu,j} f)\circ G]\ \nabla G_j +
\sum_{j \neq \mu}[(\overline \partial_{\mu,j} f \circ G)] \
\nabla \overline G_j=
\overline {\partial_\mu}\rho \ \nabla F,$$
where $\partial_{\mu,j}=\overline{\partial_\mu} \rho \ \partial_j -
\partial_j \rho \ \overline{\partial_\mu}.$

The columns (gradients) are vectors with $2n$ components.
 Solve by the Cramer's rule the resulting
$2n \times (2n-1)$ system for the $2n-1$ unknowns 
$(\overline \partial_{\mu,\nu} f) \circ G, \ \nu \neq \mu$:
\begin{equation} [\overline \partial_{\mu,\nu} f \circ G] \ \det
[\nabla G_,\nabla \overline G [\mu]]= \overline \partial_\mu \rho \ 
\det [\nabla G_1 \cdots \nabla G_n,
\nabla \overline G_1 \cdots[\mu] \cdots \nabla F \cdots
\nabla \overline G_n],
\end{equation}
where 
$\nabla F$ is inserted in the place of
the $\nu-th$ column $\nabla \overline G_\nu$ and $[\mu]$
means that the $\mu-th$ column is skipped.

Let us remind that we agreed to denote by the symbol ``$det$'' the minor
obtained by deleting the $4-th$ row corresponding to $t_3.$).
Denote $J_{\mu,\nu}$ the determinant in the right hand side of (10):
$$J_{\mu,\nu}=  \det \ [\ \nabla G; \ \nabla \overline G_1 \cdots [\mu]
\cdots \ \nabla F  \cdots \nabla \overline G_n \ ],$$
Then on $b \Sigma$ formula (10) implies:
$$J_{\mu,\nu} = K_\mu \ (\overline {\partial}_{\mu,\nu} f \circ G),$$
where
$$K_\mu= (1/\overline {\partial}_\mu \rho)
\det \ [\nabla G, \ \nabla \overline G [\mu] \ ].$$
The function $J_{\mu,\nu}$ vanishes on the central cycle $C=\{0\} \times M$
because the
first line of the determinant is formed
by the derivatives
$\partial_{\psi} F, \ \partial_{\psi} G_j, \ \partial_{\psi} \overline G_j,$
and all these derivatives vanish at $\zeta=0.$

\medskip

Now suppose that some function $J=J_{\mu,\nu}$ is not identicaly zero function
and check the properties 1, 2 of the Lemma 3.5.

First of all,  $J_{\mu,\nu}(0,t)=0$
because the first line of the determinants
consists of the derivatives in $\psi$ which vanish at $\zeta=0,$
as $\partial_{\psi}=i \zeta \partial_{\zeta}$ on holomorphic functions.

To prove that the ``phase'' part $J/\overline J$ takes the same vaues at points on $b\Sigma$ from the same $G$-fibers, we will  check that the function $K_{\mu}$ takes  only real values.
To this end, let us write, taking into account that
$\overline {\partial_\mu \rho}=\overline\partial_\mu \rho$ and
$\overline {\nabla G}=\nabla \overline G$:
$$\overline {K_\mu}= (1/\partial_\mu \rho) \ \det [\ \nabla \overline G,
\nabla G [\mu] \ ].$$
Now express the column
$\overline \nabla G_\mu$ from the relation (9) and substitute
to the above determinant. The determinants with
equal columns vanish and we obtain:
$$\overline {K_\mu}=(1/\partial_\mu \rho)
 \det
[\nabla \overline G_1 \cdots (-\partial_\mu \rho /
\overline \partial_\mu \rho) \nabla G_k \cdots
\nabla \overline G_n; \  \nabla G[\mu]].$$
Cancel $\partial_\mu \rho$ and rearrange the determinant by permutations of
the columns:
$$\overline{K_\mu}=(1/ \overline \partial_\mu \rho) (-1)^{n-1}
\det \ [\nabla \overline G[\mu];
\ \nabla G] =(1/ \overline \partial_\mu \rho) (-1)^q
\det
[\nabla G, \nabla \overline G[\mu] \ ]=(-1)^q K_\mu,$$
$q=(n-1)(n+2)/2.$

We are interested in the case $n=2$ and therefore $q=2.$
Then $\overline K_\mu=K_\mu$ and therefore
$K_\mu$ is real valued function.

For $n=2,$ the $CR-$ dimension of $\Omega$ equals 1 and the
basic tangential $\overline \partial -$ operator is,
$$\overline \partial_b=\overline \partial_{1,2}=
\overline \partial_1 \rho \ \overline \partial_2 -
\overline \partial_2 \rho \ \overline \partial_1.$$

The functions $J_{\mu,\nu},$ corresponding to the indices
$(\mu,\nu)=(1,2)$ and $(\mu,\nu)=(2,1),$ coincide and are equal to
$$J_{1,2}=J_{2,1}=J=\det \ [\nabla G_1, \  \nabla G_2; 
\ \ \nabla F \ ].$$

As in the previous section, we define
$$\Theta = \frac{J}{\overline J}.$$
We have
\begin{equation}
J=\Bigg \{ \begin{array}{cc} 
K_1 (\overline {\partial}_{1,2}f  \circ G), \ \ \overline {\partial}_1 \rho 
\neq 0  \\ 
K_2 (\overline {\partial}_{2,1}f  \circ G), \ \ \overline {\partial}_2 \rho 
\neq 0
\end{array}, 
\end{equation}
on  $b\Sigma.$  We have shown that the functions
$K_1$ and $K_2$ are real, hence 
on the manifold $b \Sigma$ holds
$$\Theta= \sigma \circ G,$$
where $$\sigma=-\frac{\overline \partial_{1,2}f}
{\overline{\overline \partial_{1,2}}f}=-\frac{\overline \partial_{2,1}f}
{\overline{\overline \partial_{2,1}}f}.$$
We have used here that 
$\overline \partial_{1,2} f=-\overline \partial_{2,1}f.$
Let us note also that the function $J(\zeta,t)=J_{1,2}(\zeta,t)$ 
is analytic in
$\zeta,$ as well as the function $J$ constructed in the previous section.

The representation (11) and definitions of the functions $K_1$ and $K_2$
imply that the set 
$J^{-1}(0)\cap b\Sigma$  contains the critical set $Crit(G)$ of $G$
on $b \Sigma.$ 

Thus, the constructed function $J$ has all the properties from Lemma 3.5.
Therefore Lemma 3.3 applies to the function $J,$ after checking the structure
of the set $J^{-1}(0)$ (condition 1 in Lemma 3.3.)


\bigskip
\noindent
{\bf 3.7. The structure of the set $J^{-1}(0)$ and
the regularization of the integrals.}

\medskip
In the previous we have constructed the Jacobians $J$
 for 
both cases $n=1$ and $n=2$
and shown that in each case the function $J$ satisfies Lemma 3.3. 
However Lemma 3.3
assumes that the zero set $J^{-1}(0)$ is a $(2n-1)$-chain for 
integrals of corresponding differential forms
over $J^{-1}(0)$ be defined.  

Thus, to apply the symmetry relation proved in Lemma 3.3 we have to understand the structure of the zero set
of the functions $J(\zeta,t).$
This is essentially the only point where
we use the assumption about real-analyticity of the mapping $G$ and
of the function $f$ (and therefore of its analytic extension, $F $).

In this case the Jacobian $J(\zeta,t)$ is real-analytic 
in $\overline \Delta \times M,$ holomorphic in $\zeta$ in the unit disc
$\Delta$, and the structure of the zero set $J^{-1}(0) $  
is well understood. 

Namely, for each value of the parameter $t$ the function
$J(\zeta,t)$ has finite number of zeros $\zeta=\zeta_j(t)$, unless 
$J(\cdot,t) \equiv 0.$ 
The zeros
constitute a collection of
points in the disc $\Delta_t =\Delta \times M$ each of which
changes continuously till it reaches the boundary circle $|\zeta|=1.$

If for some $t_0$ the zeros are not isolated then $J(\zeta,t_0)=0$ 
identically.
In this case we will call the disc $\Delta_{t_0}$  {\it zero-disc.}
Nevertheless, by real-analyticity, the zero discs are isolated.

Thus, the zero set $J^{-1}(0)$ consists of
finite number of smooth $(2n-1)-$ manifolds
$C_j=\{(\zeta(t),t): t \in M \}, \ \ \partial C_j \subset b \Sigma,$
union with and a collection (maybe empty)  of zero discs:
$$J^{-1}(0)= (\cup_j C_j) \cup Z(T),$$
where
$$Z(T)= \cup_{t \in T} \overline \Delta_t= \overline \Delta \times T.$$

For $t \notin T$ the section $J^{-1}(0) \cap \{t=t_0\}=
\{\zeta_1(t),\cdots,zeta_j(t) \}$ is finite and the function $J(\cdot,t)$
has a constant multiplicity $\kappa_j$ at its zero $\zeta_j(t).$

The set $T,$ parametrizing the zero discs, can be defined as
$$T=\cap_{s=1}^{\infty} \{ t \in M: \partial_{\zeta}^sJ(0,t)=0 \}$$
and is an analytic set.
We want to prove that the integral over the part $Z(T)$ of the critical set
$J^{-1}(0)$
contributes nothing in the integral over $J^{-1}(0)$ in formula (2).

\bigskip
\noindent
{\it Removing the zero discs in the case A, (2n,k)=(2,2).}

Consider the case $n=1.$ Then $M$ is a curve, the parameter 
$t$ is one-dimensional, and
since $J(\zeta,t)$ is a real-analytic nonzero function , it can vanish
on at most finite number of discs, that is $\# \ T < \infty.$
Let $T=\{t_1,\cdots,t_N\},$ where by $t_i$  we understand (real) coordinates
of the corresponding points in $T.$
By real-analyticity, $J$ can be represented as
$$J(\zeta,t)=q(t) J_0(\zeta,t).$$
Here $q(t)=(t-t_1)^{k_1} \cdots (t-t_N)^{k_N}$ and
the function $J_0(\zeta, t)$ possesses all the properties of the function $J$
except that the zero set $J_0^{-1}(0)$  contains no disc $\Delta_t.$

Since the parameters $t_1, \cdots t_N$   are real, we have
$$\Theta=\frac{J}{\overline J}= \frac {J_0}{\overline J_0}.$$
where $J_0$ has on each disc $\overline \Delta_t$ only finite number of 
zeros of finite order. Then $J_0$ satisfies the conditions of Lemma 3.3
and formula (2) holds for $J_0.$ Thus, in the presence of zero disc
we compute integrals in formula (2) for the ``regular'' part, $J_0$ of $J.$


\medskip
\noindent
{\it Removing the zero discs in the case B, (2n,k)=(4,3).}

In this case $\dim M=3$ and the dimension of the analytic
subset $T \subset \Omega$ can be 0, 1 or 2. By the (real) dimension
we understand the maximal dimensions of Whitney strata (see \cite{GM}.)

If $\dim T=0$ then $\dim Z(T)=2 .$ Since in (2) we integrate
3-form, the cycle  $G(Z(T))$ is negligible.

Suppose $\dim T=1.$ Then $Z(T)$ contains 3-chains
of the form $c= \overline \Delta \times \gamma,$ where
$\dim \gamma=1.$  This chains might contribute in (2).

However, the condition $H_1(M)=0$
in Theorem 2.2 implies that $\gamma$
is homological (relatively homological) to zero in $M.$

Then the cycle $G(c)$ is homological to zero in
the image $G(\overline \Sigma)$ and since the Martinelli-Bochner form is
closed, the integral over the cycle $G(c)$ in (2) is zero.

Finally, if $\dim T=2$ and $\gamma \subset T$ is a stratum of
the pure dimension 2,
then locally, in a neighborhood of a point $t_0 \in M,$ we have
$$\gamma=\{q(t)=0\},$$
where $q(t)$ is real-valued. Then, as in the case $n=1$,
locally $$J(\zeta,t)=q(t)J_0(\zeta,t)$$
where $J_0(\zeta,t)$ has only isolated zeros on discs
$\overline \Delta_t$.

Since $J/\overline J=J_0/\overline J_0,$ then we conclude that
the set $c=\overline \Delta \times \gamma$ is removable ,
by the same argument as for $n=1.$

\section {End of the proof of Theorem 2.2
for the case $\partial M=\emptyset.$}

Everywhere in this section $n=1$ or $n=2.$
\bigskip
Our strategy is as follows.
We assume that the mapping $Q$
is nondenerate. Then the Jacobi minors $J$ constructed
in Section 3.5 and
3.6 are not identically zero.
Using the symmetry relation in Lemma 3.3,
for the case when the curve $G^{-1}(b)$ is just empty,we obtain
that the integral of the differential form $(G-b)^* \beta_{MB}$
over the set
$J^{-1}(0)$ is zero. On the other hand, this set contains the
nontrivial cycle
$C$ which is mapped by $G$ to a nontrivial cycle in the image. Hence the
integral can be made nonzero by appropriate
choice of the point
$b$ and we obtain contradiction with the 
assumption that $Q$ is nondegenerate.
\bigskip
Let us start realizing this program.
Let
 $J$ be the function constructed in Sections 3.5 and 3.6 and
assume that $J \neq 0.$
Now we want to apply Lemma 3.3 to the function $J.$

All the conditions
of Lemma 3.3 are fulfilled. Removing the zero discs according to
the regularizations procedure in Section 3.7, we represent
the $(2n-1)$-chain $J^{-1}(0)$ as the union of connected smooth
manifolds with multiplicities:
$$J^{-1}(0)=\cup_j \kappa_j C_j.$$
One of these manifolds is just the cycle $C^{\prime}=C_{j_0}.$
The multiplicities $\kappa_j$ are the winding numbers (indices) of
the function $J$ at its zeros with respect to $\zeta$ (see Section 3.7).

Let $b \in \mathbb C^n$ is such that
$J^{-1}(0) \cap G^{-1}(b) =\emptyset.$
Apply Lemma 3.3:
$$2 \int_{J^{-1}(0)}(G-b)^*\beta_{MB} =
\frac{1}{2\pi i} \int_{G^{-1}(b)} (\frac{dJ}{J}-
\frac{d\overline J}{\overline J}).$$

After change of variables $(z_1,z^{\prime})=(J(\zeta,t),
G(\zeta,t)-b),$ the latter identity transforms to
\begin{equation} 2 \sum_j \kappa_j \int_{G(C_j)} \beta_{MB}(z^{\prime} -b)=
\frac{1}{2\pi i} \int_{G^{-1}(b)} (\frac{dJ}{J}-
\frac{d \overline J}{\overline J}).
\end{equation}

\begin{Lemma}
In the case  $(2n,k)=(2,2)$ and $(2n,k)=(4,3),$
the numbers $\kappa_j$  are positive.
\end{Lemma}
\begin{proof}
Each chain  $C_j \subset J^{-1}(0)$ is defined, near generic points,
by the equation
$\zeta=\zeta_j(t), \ t \in M.$

In Sections 3.5 and 3.6 we have constructed
the functions $J=J_-(\zeta,t)$ and $J=J_{1,2},$
for the cases $(2,2)$ and $(4,3)$ correspondingly.
In both cases these functions are analytic in $\zeta,$
for the determinants defining $J$ do not contain antiholomorphic columns.

Hence the winding numbers $\kappa_j= \kappa_j (t)$ of the
mapping $\zeta \mapsto J(\zeta,t)$
on a small circle \ $|\zeta - \zeta(t)|=\varepsilon,$ \
at an isolated zeros
$\zeta=\zeta_j(t), $ are positive integers, equal to
the muplitipicities of the zero of the holomorphic function $J$.
If $\zeta_j(t)$ consists of the boundary zeros,
$|\zeta_j(t)|=1$ then $\kappa_j$ must be taken 1/2 of the multiplicity of the
zero $\zeta_j(t)$ because the integration in (11) is performed within
the unit disc.

\end{proof}
\begin{Remark}
Lemma 4.1 is the only point where we use the restriction
$n \le 2$ for the dimension. This restriction means that in both
cases, A and B, of Theorem 2.2, the $CR-$ dimension of the manifold $\Omega$
is one. In this case, the Jacobians $J$ in Sections 3.5 and 3.6
are sense-preserving mappings
and the winding numbers $\kappa_j$ are positive,
which is an important point in the proof.
\end{Remark}

\bigskip
Now we can finish the proof.
We know that $J$ vanishes on a homologically nontrivial cycle
$C$ in $\overline \Sigma:$
$$ C \subset J^{-1}(0),$$
so that $C$ is one of the cycles $C_j$ constituting the chain
$J^{-1}(0).$

Let $b \notin G(\overline \Sigma)$.
Apply the identity (12) to that value $b$.
Since $b$ in not in the image of $G$ then
$G^{-1}(b)=\emptyset$ and hence the integral
over the curve $G^{-1}(b)$ is zero:
\begin{equation} \int_{G^{-1}(b)} (\frac{dJ}{J}-
\frac{d \overline J}{\overline J}) =0.
\end{equation}
On the other hand, according to Lemma 3.6,
the image cycles $G(C_j)$ are all cooriented and
the multiplicities $\kappa_j$ are all positive,
$ \kappa_j > 0.$

Moreover, by condition (iii) of Theorem 2.2,
the cycle $c=G(C)$ is not homological to zero in
$G(\overline \Sigma).$ However it is homological to zero in $\mathbb R^{2n}$
and therefore $c=\partial c^{\prime}$ for some $2n-$ dimensional cycle
$c^{\prime}$ in $\mathbb R^{2n}.$  This cycle does not belong entirely to
$G(\overline \Sigma)$ and therefore  the point $b \in \mathbb C^n \setminus
G(\overline \Sigma)$ can be chosen
so that $c \neq 0$ in $H_{2n-1}(R^{2n} \setminus \{b\}).$

We have the isomorphisms of the homology groups with coefficients 
in $\mathbb R:$
$$H_{2n-1}(\overline \Sigma) \cong H_{2n-1}(M) \cong \mathbb R \cong
H_{2n-1}(R^{2n} \setminus \{b\}).$$

Martinelli-Bochner form $\beta_{MB}(z^{\prime}- b) $ is a closed nonexact form
in $\mathbb R^{2n} \setminus \{b\},$  
representing the generator  in 
de Rham cohomology group $H^{2n-1}(\mathbb R^{2n} \setminus \{b\}, \mathbb R)
\cong \mathbb R.$
By the duality,
the integral $$\int_{G(C_j)} \beta_{MB}(z^{\prime} - b)=m_j,$$
where $$m_j = \{G(C_j)\} \in H_{2n-1}(\mathbb R^{2n} \setminus \{b\})$$ 
is the homology class of the cycle $G(C_j).$
Since the cycles $G(C_j)$ all are cooriented, we have $m_j \ge 0$ and
$m_j=0$ if and only if the cycle $G(C_j)$ is homological to zero in
$\mathbb R^{2n} \setminus \{b\}$.

The cycle $C$ is one of the $C_j's,$ i.e. $C=C_{j_0}.$  
Then $m_{j_0} > 0$ because, as we saw,
$c=G(C)$ is a homologically nontrivial cycle in
$\mathbb R^{2n} \setminus \{b\}$  and
hence the total sum
$$
\sum_j \kappa_j \int_{(G-b)(C_j)} \beta_{MB}
\ge \kappa_{j_0} m_{j_0} >0.$$

Together with (13) this inequality contradicts to
(12).

We have assumed that $J \neq 0$ and arrived
to contradiction.Therefore  $J \equiv 0.$ Now are ready to
obtain the final conclusion of Theorem 2.2.

\medskip
\noindent
{\bf The final step: proof that $Q$ is degenerate, $F$ is constant of $G$-fibers 
and $f$ is a $CR$ function.}

The contradiction we just have obtained says that the minor
$J$ (the case A,  n=1) or the minors $J^0,J^1,J^2,J^3$ (the case B, n=2) 
of the Jacobi matrix $J(F,G)$ with respect to variables $\psi,t$
vanish identically in our main manifold  $\overline \Sigma$. 
Observe that since the functions $F,G$ are analytic in $\zeta$,
the minors built from the derivatives in $\zeta$  vanish as well.

In the case $n=1$, when the parameter $t$ is one-dimensional,   
we obtain that $rank_{\mathbb C} J(F,G) < 2$ and $Q$ is degenerate
in accordance with the claim of Theorem 2.2. 

Moreover, $J=0$ implies that
$F=\hat f \circ G$ in $\Sigma$ for some  function $\hat f$. Indeed,
we have shown in Section 3.5 that
the Jacobian $J$ is proportional to the directional derivative of $F$
along the level curves $G=const$ and hence $J=0$ means that
$F=const$ on $G$-fibers.

Also the the relation between $J$
and the $\overline {\partial}$- derivatives of the function $f$ proved in
Section 3.5
implies 
$$\overline \partial f \circ G =0$$ 
on $b\Sigma$
which means that $f$ is holomorphic in $\Omega=G(b\Sigma).$

\bigskip
Now turn to the case $n=2.$ Again, vanishing    the minors 
$J^k, k=0,\cdots,3,$
implies that $rank_{\mathbb C} J(F,G_1,G_2) < 3$ and hence $Q$ is degenerate.
Moreover, since due to the regularity condition, $\nabla G_1,\nabla G_2$
are linearly independent over $\mathbb C$ at any point in $\Sigma$,
we have $\nabla F \in span\{\nabla G_1,\nabla G_2 \}.$ If $\zeta=\zeta(s),
t_k=t_k(s), k=1,2,3; s$
is one-dimensional parameter,
 is a  level curve  $\Gamma=G^{-1}(b)$, then the directional derivative along the level curve is
$$
\partial^{G}F=\zeta^{\prime}(s)\partial{\zeta F}+
t_1
^{\prime}(s)
\partial_{t_1}F + t_2^{\prime}(s) \partial_{t_2}F + 
t_3^{\prime}(s) \partial_{t_3}F.$$
Since $\partial^{G}G_1=\partial^{G}G_2=0$, we conclude from the
above linear dependence of the gradients:
$$\partial^{G}F=0$$
and therefore $F=const$ on any $G$-fiber.

Finally,    the     identity (11) implies  that for tangential 
$\overline{\partial}$-derivative holds 
$$\overline{\partial}_b f \circ G=0$$ 
on $b\Sigma$, and this means that
$f$ is a $CR$-function on $\Omega.$

\bigskip
\noindent
\section{End of the proof of Theorem 2.2 for
the case $\partial M \neq \emptyset.$}

The proof for the case of nonclosed parameterizing manifold $M$
follows the same line
as in the case of $M$ with the empty  boundary. 

Again, we assume that the Jacobian $J$ in the case $n=1$ or one of the
three minors $J^1,J^2,J^3$ in the case $n=2$ is not identical zero.
In all cases, we denote this determinant by $J$. Our goal is
to arrive to a contradiction which would allow to conclude that
all the above minors vanish identically. Then all the conclusions of Theorem 2.2 follow by the same arguments as in the previous Section.

The main difference with the case of the previous Section,
where the parametrizing manifold $M$ was closed,
is that now we have to deal with
the relative homology groups. 

Our basic manifold now is also a bit different. Namely, we set
$$Y = G(\overline \Sigma) \setminus G(\overline \Sigma_0), \ 
 \Sigma_0=\Delta \times \partial M $$ and define
$$\Sigma{^\prime}= G^{-1}(Y).$$

Thus, we take only that part of $\Delta \times M$ which is
the full $G-$ preimage of its image.
Note, that by the condition the set $Y$ is nonempty and connected.
The topological boundary of $\Sigma^{\prime}$ is
$$\partial \Sigma^{\prime}= (\partial \Sigma \cap G^{-1}
(\overline Y)) \cup  G^{-1}(\partial Y).$$
Taking into account that, in turn,
 $\partial \Sigma= (\partial \Delta \times M) \cup
(\Delta \times \partial M)=b \Sigma \cup \Sigma_0$ and
that $G^{-1}(Y) \cap \Sigma_0=\emptyset,$
 we obtain that
the boundary of the new manifold is
$$\partial \Sigma^{\prime} =
(G^{-1}(Y) \cap b \Sigma) \cup  (G^{-1}(\partial Y)\cap   \Sigma_0) .$$

Thus, we assume that the mapping $Q=(F,G)$ is not degenerate and,
as we have mentioned already, 
construct the Jacobian minor $J$ and the function $\Theta=J/\overline J$, 
as in Sections 3.5 and 3.6. Our  main assumption about $Q$
is that the function $J$ is not identically zero.

Assume initially that the zero set $J^{-1}(0)$ does not contain
entire discs $\Delta_t.$ Proceed as in the case
$\partial M=\emptyset$ (Lemma 3.3). Namely, we delete from $\Sigma$
small neighborhoods of the zero
sets and apply Stokes formula to the remainder. Again, we show as in Lemma 3.3,
that the surface
integral contributes nothing, because of degeneracy
of the mapping
$(J,G)$ on $\partial \Sigma^{\prime}.$
However, an
extra surface term appear,
corresponding to the extra part of the boundary. As the result, we obtain:
\begin{equation}  2 \int_{J^{-1}(0)} (G-b)^*\beta_{MB} =
\frac{1}{2\pi i} \int_{G^{-1}(b)}(\frac{dJ}{J}-\frac{d \overline J}
{\overline J})
 - \frac{1}{2\pi i} \int_{G^{-1}(E)}
(\frac{dJ}{J} - \frac{d \overline J}{\overline J}) \wedge (G-b)^*\beta_{MB},
\end{equation}
where
$$E=\partial Y \cap G(\Sigma_0)$$
and correspondingly $G^{-1}(E)$ is
the $2n-$
dimensional extra part of the boundary of $\Sigma^{\prime}$, coming from
$\partial M.$

Note that the set $J^{-1}(0)$ is not necessary closed cycles any longer, 
but it is a cycle relatively to
$G(\Sigma_0).$  Then we will use duality arguments based on counting of 
intersection indices, in place
of computing periods by integration of Martinelli-Bochner type 
integrals as in the case $\partial M=\emptyset.$

Write (14) as
$$\chi(b)=Z(b)\ + \ N(b),$$
where $\chi(b)$ stands for the left hand side in (14) and
$$Z(b)=\frac{1}{2\pi i} \int_{G^{-1}(b)}(\frac{dJ}{J}-
\frac{d \overline J}{\overline J}), \ \
N(b)= - \frac{1}{2\pi i} \int_{G^{-1}(E)}
(\frac{dJ}{J}-\frac{d \overline J}{\overline J})\wedge (G-b)^*\beta_{MB}.$$

By change of variables,
$$\chi(b) =\frac{1}{\pi i} \sum_j \kappa_j \int_{G(C_j)}\beta_{MB}
(z^{\prime} -b) ,$$
where $C_j$ are connected components of the $(2n-1)$-chain $J^{-1}(0)$
(as in  the case $\partial M=\emptyset$),
and $\kappa_j$
are their multiplicities with respect to $J.$ According to Lemma
4.1 the numbers
$\kappa_j$ are positive.
\begin{Lemma} The function $Z(b)$ is integer valued as long as
$J^{-1}(0) \cap G^{-1}(b) = \emptyset.$
\end{Lemma}
\begin{proof}
By the construction in Lemma 3.5 (property 2), the function
$\Theta=J/\overline J$ is
$G-$ compatible on $b \Sigma=\partial \Delta \times M$,
i.e. it takes same values at the end points of the curves
$G^{-1}(b).$
Then the integral of the logarithmic derivative equals to the variation
of the argument of $\Theta$ along the curve $G^{-1}(b)$ and is integer:
$$Z(b)= \frac{1}{2 \pi i} \int_{G^{-1}(b)} \frac{d\Theta}{\Theta}=
(1/2 \pi) Var_{G^{-1}(b)} arg \ \Theta \in \mathbb Z.$$
\end{proof}
\begin{Lemma} The function $N(b)$ is continuous in the domain $Y$.
\end{Lemma}
\begin{proof} By the construction the sets
 $G^{-1}(Y)$ and $G^{-1}(E)$ are disjoint. Therefore
when $b \in Y$ then the
 differential  form $(G-b)^*\beta_{MB}$ has no singularities on
the surface of integration, $G^{-1}(E)$.
The 1-form
$dJ/J - d\overline J/\overline J $ has  removable singularities
at isolated zeros of $J$ on the $2n$-surface $G^{-1}(E)$ and does not
depend on $b$. It follows that the integral defining the function $N(b)$
depends continuously on $b \in Y.$
\end{proof}

Let $J^{-1}(0) =\cup_j C_j$ be the decomposition into connected chains.
We know that the set $J^{-1}(0)$ contains the central cycle
$C_{j_0}= C =\{0\} \times M,$ for which
the condition (iii) yields  that the homology class
$[G(C)] \in H_{2n-1}(G(\Sigma),G(\Sigma_0))$ is not 0.

The portion of $G(\Sigma_0)$ in $Y$
is just the set $E$ and the homological nontriviality of $G(C)$
means  that $\partial G(C) \subset E$ and
$G(C)$ can not be made the boundary of a $2n$- cycle in $Y \cup E$
by adding a $(2n-1)$-chain contained in $G(\Sigma_0).$ In other words,
the $(2n-1)$-chain $G(C)$ has its boundary in the set $E,$
but is not homologous to any $(2n-1)$-chain in $E$. The chain $G(C)$
intersects each domain $D_t=G(\Delta_t), \ t \in M.$

\begin{Lemma} There exists a
connected continuous curve, $L \subset \mathbb C^n \setminus 
\overline {G(\Sigma_0)}$,
with end points $p, \ q \in \mathbb C^n \setminus \overline{G(\Sigma)}$
intersecting the $(2n-1)$-chain $G(C)$ with the 
intersection index $m \neq  0.$ 
\end{Lemma}
\begin{proof} Denote for brevity $X=G(\Sigma), \ A=G(\Sigma_0).$
The 
chain $G(C)$ represents a nonzero element  
$$ 0 \neq [G(C)] \in H_{2n-1}(X, \overline A).$$
By Poincare-Lefshetz duality 
(see e.g. \cite{Sp}, Ch. 6), 
there exists an element 
$$ 0 \neq [L] \in H_1
(\mathbb C^n \setminus \overline A, \mathbb C^n \setminus \overline X)$$ 
with nonzero intersection index, $$ind \ (L \cap G(C)) \neq 0.$$ 
The 1-chain $L$ has the required properties.
\end{proof}

\noindent
{\bf The final step: proof that $J \equiv 0.$}

The constructed  1-chain $L$ enters and leaves the set $Y,$ avoiding the set 
$E \subset G(\Sigma)$ and crossing the chain $G(C)$
with a nonzero intersection index.
Reversing, if needed, orientation of the  path $L \subset Y$ we can assume
that the intersection index $m=m_{j_0}>0$ .

By Sokhotsky-Plemelj theorem for the Martinelli-Bochner type integrals
(see,e.g. \cite{LC}),
each time  $b \in L$ intersects transversally some chain $G(C_j),$
the integral $\chi(b)$
changes for +1 or -1, depending
on the index of the intersection. 
If $m_j$ is the total variation of $\chi(b),$
resulted from the crossing the cycle $G(C_j),$ 
then numbers $m_j$ are of the same sign since
all the cycles $G(C_j)$ are cooriented.

Thus, the algebraic number  of jumps 
the function $\chi(b), \ b \in L,$ makes, after passing through
$Y,$ is at least $\kappa_{j_0}m_{j_0}$:
$$\#jumps_{b \in L} \chi(b) \ge \kappa_{j_0}m_{j_0}
>0.$$
But $\chi(b)=N(b) + Z(b), \ b \in Y$ and  the function
$N(b)$ changes continuously. On the other hand, by Lemma 5.1,
 the function $Z(b)$
is    integer valued. Therefore $Z(b)$
changes  only when  $b$ crosses the set $J^{-1}(0)$
which results in a jump of
value of $\chi(b)$. Thus, the total algebraic number of jumps
 of $\chi(b)$ and $Z(b)$ are
equal and hence:
$$\#jumps_{b \in L} \ Z(b)=\# jumps_{b \in L} \chi(b) 
\ge \kappa_{j_0}m_{j_0}>0.$$

However, $G^{-1}(b)=\emptyset$ when $b \notin Y$. Therefore
the function $Z(b)$ takes the value 0 before the argument $b$
enters the set $Y$  and after it has left $Y$. In
other words,
if $b_0,\  b_1 \notin Y$ are the starting and end points of the path
$L$ then  $Z(b_0)=Z(b_1)=0$ and therefore
$$Var_{b \in L} \ Z(b) =0.$$
This contradiction says that $J=0.$

It remains to explain how to get rid of the zero discs.
For the case $(2n,k)=(2,2)$ the argument from  Section 3.7
works as it does not depend on whether $M$ is closed or not.

Let $(2n,k)=(4,3)$ and let $T$ and $Z(T)$ are sets defined in Section 3.7.
When              $\dim T=0$     and            $\dim T=2,$
the zero discs are removable for the same reasons as in Section 3.7.
When $\dim T=1$ and $Z(T)$ contains 3-chain of the form $
c=   \overline \Delta \times \gamma,$
then by the condition $H_1(M,\partial M)=0$
the relative 1-cycle $\gamma$ is relatively homologically zero. This implies
that
the total intersection index  of the 1-chain $L$ with
the relative 3-cycle $G(c)$ is zero and therefore the cycle $G(c)$
contributes nothing in the result of our computation.

The proof of Theorem 2.2 is completed.

\section{Proof of Theorem 4.}


We restrict ourselves by the proof for 
the case $\partial M=\emptyset,$ i.e. $M=S^1.$
Let $Q(\zeta,t)=F(\zeta,t), G(\zeta,t), \ (\zeta,t) \in \Sigma=
\Delta \times S^1,$ be a 
parametrization of the family  $D_t$ in the formulation of Theorem 4. 
We have $Q(b \Sigma)=\Lambda$ because the curves 
$\partial D_t, \ t \in S^1,$ cover $\Lambda.$ 
By regularity, the mapping 
$$Q: b \Sigma \mapsto \Lambda$$
is finitely sheeted covering by the torus 
$b \Sigma =T^2$ of $\Lambda \setminus \partial \Lambda.$ 
The image $Q(\Sigma)$ is the union of the discs
$D_t$ and its dimension is at most 3.

We want to prove that the mapping $Q=(F,G)$ is 
degenerate in $\Sigma,$ and correspondingly,
$Q(\Sigma)$ is  2-dimensional and  entirely contained in  $\Lambda.$

We  do so by slight modification of the arguments in the proof 
of Theorem 2.2. 
First of all we need to prove the symmetry relation analogous to (2).
 
To this end we will check first that, as in Lemma 3.5,  the Jacobian 
$$J= \det [ \nabla F, \ \nabla G]= \frac{\partial(F,G)}{\partial(\psi,t)},$$
with respect to local coordinates $\psi, t$ on $b \Sigma$ has constant 
"phase"  on the $Q-$fibers 
$Q^{-1}(b) \cap b\Sigma, \ b \in \mathbb C^2.$ 
More precisely, we want to prove that $J/\overline J$ is constant 
on the above fibers, outside of the zeros of $J:$

\begin{Lemma} Let $u_1, \ u_2 \in b \Sigma$ belong to the same 
$Q$-fiber: \ $Q(u_1)=Q(u_2).$
If $J(u_1) \neq 0$ then $J(u_2) \neq 0$ and
$$\frac{J(u_1)}{\overline{J(u_1)}}=\frac{J(u_2)}{\overline {J(u_2)}}.$$
Therefore on $b \Sigma \setminus J^{-1}(0)$ holds
$$\frac{J}{\overline J}=\sigma \circ Q,$$
for some function $\sigma.$ 
\end{Lemma}
\begin{proof} Each pair of vectors in $\mathbb C^2$
$$\{\partial_{\psi}Q(u_j), \partial _t Q(u_j)\}, \ j=1,2,$$
belongs  to the tangent space $T\Lambda_b$ at the point $b=Q(u_1)=Q(u_2).$ 
Since the Jacobian $J$ does not vanish at the points $u_1, u_2$, each pair 
constitutes a basis in the 2-dimensional
real space $T\Lambda_b.$ Let $\mathcal {A}$ be the (real) transition matrix 
from one basis to another.   
The Jacobians are related by
$$ \frac{\partial(F,G)}{\partial(\psi,t)}(u_1)=\det \mathcal {A} \cdot 
\frac{\partial (F,G)}{\partial (\psi,t)}(u_2),$$
that is $J(u_1)=\det \mathcal {A} \cdot J(u_2)$ and lemma follows as 
$\det \mathcal {A}$ is real and nonzero and hence it cancels
when one divides $J$ by $\overline J.$
\end{proof}

Now we need an analog of Lemma 3.1.
\begin{Lemma} The Brouwer degree of the mapping $Q: T^2= b \Sigma \mapsto 
\Lambda \subset \mathbb C^2$ equals zero.
\end{Lemma}
\begin{proof} Lemma 6.2 folows from the condition 
$\partial \Lambda \neq \emptyset$
and from the fact that any smooth mapping of a compact closed 
manifold 
to a compact manifold with nonempty boundary has the Brower degree 0. 
Indeed, the image
$Q(T^2)$ can be slightly homotopically contracted inside $\Lambda$ so that some
points in $\Lambda$ would have empty preimage with respect to a mapping 
$Q^{\prime}$,
homotopically equivalent ot $Q.$ Therefore, by local definition of the degree,
$\deg Q= \deg Q^{\prime}=0.$ 
\end{proof}

Now assume that $J$ is not identically zero and correspondingly 
$J^{-1}(0)$ is a 1-chain (we remove the zero discs exactly as
we did in Section 3.7).
We need the analog of Lemma 3.3 in the form
\begin{Lemma} Let $\omega$ be a closed differential 
1-form in a neighborhood of $Q(\overline \Sigma).$
Then 
\begin{equation}
\int_{J^{-1}(0)} Q^* \omega =0.
\end{equation}
\end{Lemma}
\begin{proof}  We follow the proof of Lemma 3.3.
Define $$\Theta=\frac{J}{\overline J}.$$
The form $$\Xi=\frac{d\Theta}{\Theta} \wedge Q^*\omega$$
is closed in $\overline \Sigma \setminus J^{-1}(0).$
The differential form $$\frac{d\Theta}{\Theta}=\frac{ dJ}{J} - 
\frac{d\overline J}{\overline J}$$ computes
the current $d\Theta/\Theta=
2[J^{-1}(0)]$. Lemma 6.1 and Lemma 6.2 imply 
that $\Xi$ integrates to zero on $b \Sigma \setminus J^{-1}(0)$
and then Stokes formula leads to  identity (15), in the same 
way as in Lemma 3.3.
\end{proof} 

Now we can complete the proof of Theorem 4. Since $J(0,t)=0,$ 
the zero set $J^{-1}(0)$ contains the central cycle C= $\{0\} \times S^1.$ 
The image $Q(C)$ meets each analytic disc $Q(\Delta \times S^1)$ and 
by the condition of homological nontriviality
the cycle $Q(C) \subset Q(\Sigma)$ is not homologous to zero. 
By de Rham duality, there exists a closed diferential 1-form
$\omega$ on $Q(\Sigma)$ such that
\begin{equation}
\int_{Q(C)} \omega \neq 0.
\end{equation} 
The final argument is as in Section 4. 
Namely, the mapping $J$ is holomorphic in $\zeta$ and 
hence all cycles $G(C_j),$
where $J^{-1}(0)=\cup C_j,$ are cooriented. Hence the integrals 
of $\omega$ over $Q(C_j)$ are of the same sign and since the cycle $C$
is one of $C_j's$, (15)
implies $$\int_{Q(C)} \omega =0.$$
This contradicts to (16). Therefore
$$i \zeta \ \frac{\partial(F,G)}{\partial(\zeta,t)}= 
\frac{\partial(F,G)}{\partial(\psi,t)} =0$$
identically. 

The vectors
$$\partial_{\psi}(F,G)(u), \ \partial_t (F,G) (u), \ u \in b\Sigma,$$
span the tangent space $T\Lambda_b, \ b=Q(u)=(F(u),G(u))$. 
We have shown that the Jacobian of $Q=(F,G)$ vanishes identically and hence 
the above two vectors are linearly
dependent over $\mathbb C.$ Therefore 
the real 2-dimensional tangent space $T\Lambda_b$ at each point 
$b \in \Lambda$ is  a complex line in $\mathbb C^2.$
This proves that $\Lambda$ is a 1-dimensional 
complex manifold. By the uniqueness
theorem for analytic functions, the attached analytic discs 
entirely belong to the manifold: $D_t \subset \Lambda,
t \in S^1.$ Theorem 4 is proved.

\section {Concluding remarks.}

\begin{itemize}

\item {\bf Comments on the proof.}

The proof of Theorem 2.2  for the case $\partial M \neq \emptyset$
given in Section 5 is also valid for the case $\partial M=\emptyset$.
As a matter of fact, for the case
of closed manifold $M$ both proofs coincide as  one can count the number
of nontrivial $(2n-1)$- cycles either by computing logarithmic residue
type integrals, as in Section 4, or by counting of intersection indices in
transversal direction, as in Section 5.

Some ingredients of the proof in Section 5
are close to the argument, presented in a different form
and for the case of
circles in the plane, in Tumanov's article \cite{T2}.

Let us comment on this.
In \cite{T2} (as well as in \cite {AG},\cite {T1}) the Schwarz 
function  of the circles (the meromorphic extension of $\overline z$
from the circle in the corresponding disc) is a key tool.
Globevnik observed in \cite{G7} that the construction in \cite{T2}, 
works also for families parameterized by mapping of the form
$G(\zeta,t)=g(\zeta)+it, \ t \in \mathbb R,$
where $g$ is a conformal mapping of the unit disc, having the symmetry
$g(\overline \zeta)=\overline {g(\zeta)}.$
This allowed him to obtain the corresponding test of analyticity
for a special families of non-circular curves.

In fact, the Shwarz function in the above articles plays the role similar to
that of Jacobian $J$ in our approach. To explain this,
observe that
the Schwarz function of the circle
coincides with $\overline z$ on the circle
$|z-a|=r$ and therefore is $G$-compatible for the parameterizing mapping
$G(\zeta,t)=a(t) + r(t)\zeta$. This function
develops a simple pole at the center
of the circle which can be turned to zero by taking the reciprocal function.
It results in the function satisfying the two conditions: it is 
$G$-compatible on $b\Sigma$ and vanishes on the central cycle $C$.
Existence and construction of nontrivial functions with these two properties is just
the key ingredient in our proof. 

Nevertheless, the only domains,
for which the Schwarz function
has a simple pole, are discs. In fact, there is a larger class of domains with
meromorphic Schwarz functions, so called quadrature domains (see,e.g.
\cite{Sh}), but this class is pretty restrictive, too, as even simple domains,
like ellipses, do not belong to it. So, the approach based on the Schwartz
functions does not seem to work for general domains and rather works for discs 
or some other quite special families. 

The key point in our proof is that the needed
properties ($G$-compatibility and vanishing on a nontrivial cycle in $\Sigma$)
are provided by the Jacobian determinants $J$, 
which encode the holomorphy or being $CR$ of $f.$ 
The Jacobi determinants are intrinsically produced
both by the family of curves (i.e. by the parametrizing
function $G$) and by the analytic extensions $F$ of the 
original function $f$.    
This circumstance allows to deal with arbitrary functions $G$ and
makes the proof work for general Jordan curves.

\item
{\bf Smoothness assumptions.}

Our construction involves derivatives and therefore requires
at least smoothness both
of the family of the analytic discs 
and, which is perhaps even worse, of the function $f.$

Moreover, we use integration of differential forms
over the zero sets $J^{-1}(0)$ and assume even stronger condition of
real-analyticity 
to provide nice structure of the zero sets and the integrals being defined. 

In this article we have focused on the very
constructions  and therefore 
assume all objects under consideration maximally smooth.  
However, the technology of 
{\it currents} allows integration of differential forms over ``bad'' sets 
and hence, as it was noticed in the beginning of the article, 
may allow to extend our approach  
for differentiable, say $C^{\infty},$ functions $f$. 

Another possible way of relaxing smoothness assumptions might be the following.
Observe that the function $f$ in Theorem 1 (the strip-problem) 
appeared, after Theorem 1 is proved,
real analytic in $\Omega$ because it is holomorphic there.
Therefore one could try to derive real-analyticity of $f$ 
from the original conditions, under the initial minimal assumption of
continuity of $f$. 

It seems plausible that it might be 
done by using argument from \cite {T1} based on the egde of the wedge theorem. 
This argument might lead to real analyticity of $f$ in a narrow open
curved strip in $\Omega$ which would be enough to conclude that all the
zero discs in $J^{-1}(0)$ are isolated. The latter is important
for integrals of differential forms in our proofs to be defined.
It should be mentioned,either, that the discussed method may work 
only for the strip-problem, since in higher dimensions (for instance, in Globevnik-Stout conjecture) one can not
expect {\it a posteriori} real analyticity of $f$.
We plan to return to the problem of relaxing the 
smoothness assumptions elsewhere.

\item
{\bf The case $n >2.$}

The restriction $n \le 2$ for the dimension is used only at
one, albeit crucial, point, namely for checking that the multiplicities
$\kappa_j$ are of the same sign. This is provided by
analyticity of the Jacobians
$J$ in $\zeta,$ which is true just
in the cases $n=1, \ \ \dim_{\ \mathbb R}\ \Omega=2,$ \ \
and   $n=2,\  \ \dim_{\ \mathbb R}\ \Omega =3.$
Nevertheless, the method developed in this article can be applied, 
{\it mutatis mutandis,} in higher dimensions by considering minors of 
less dimensions. Then the proof works for families of attached analytic discs
satisfying certain conditions for less dimensional subfamilies.
We are going to return to the higher dimensional case elsewhere.

\bigskip

\end{itemize}

{\bf Acknowledgments}

I am greatly indebted to Josip Globevnik for numerous stimulating
and inspiring discussions of the strip-problem. 
I thank Alexander Tumanov for informing 
me about  the   preprint of his article
\cite{T2}. I am grateful to Steve Shnider who
read the first version of the manuscript and made several  useful remarks. 

\bigskip
The main results of this article were announced in the short note \cite{AgCR}.

\vskip.15in

\end{document}